# The uncountable spectra of countable theories

By BRADD HART, EHUD HRUSHOVSKI, and MICHAEL C. LASKOWSKI *

### Abstract

Let $T$ be a complete, first-order theory in a finite or countable language having infinite models. Let $I(T, \kappa)$ be the number of isomorphism types of models of $T$ of cardinality $\kappa$. We denote by $\mu$ (respectively $\hat{\mu}$) the number of cardinals (respectively infinite cardinals) less than or equal to $\kappa$.

THEOREM. *$I(T, \kappa)$, as a function of $\kappa > \aleph_0$, is the minimum of $2^\kappa$ and one of the following functions*:

1. $2^\kappa$;

2. *the constant function* 1;

3. $\begin{cases} |\hat{\mu}^n/\sim_G| - |(\hat{\mu} - 1)^n/\sim_G| & \hat{\mu} < \omega; \quad \text{for some } 1 < n < \omega \text{ and} \\ \hat{\mu} & \hat{\mu} \geq \omega; \quad \text{some group } G \leq \mathrm{Sym}(n) \end{cases}$

4. *the constant function* $\beth_2$;

5. $\beth_{d+1}(\mu)$ *for some infinite, countable ordinal d*;

6. $\sum_{i=1}^d \Gamma(i)$ *where $d$ is an integer greater than 0 (the depth of $T$) and*

   $\Gamma(i)$ *is either* $\beth_{d-i-1}(\mu^{\hat{\mu}})$ *or* $\beth_{d-i}(\mu^{\sigma(i)} + \alpha(i))$,

   *where $\sigma(i)$ is either $1, \aleph_0$ or $\beth_1$, and $\alpha(i)$ is 0 or $\beth_2$; the first possibility for $\Gamma(i)$ can occur only when $d - i > 0$.*

The cases (2), (3) of functions taking a finite value were dealt with by Morley and Lachlan. Shelah showed (1) holds unless a certain structure theory (superstability and extensions) is valid. He also characterized (4) and (5) and

*The first author was partially supported by the NSERC. The third author was partially supported by the NSF, Grant DMS-9704364. All the authors would like to thank the MSRI for their hospitality.

showed that in all other cases, for large values of $\kappa$, the spectrum is given by $\beth_{d-1}(\mu^{<\sigma})$ for a certain $\sigma$, the "special number of dimensions."

The present paper shows, using descriptive set theoretic methods, that the continuum hypothesis holds for the special number of dimensions. Shelah's superstability technology is then used to complete the classification of the all possible uncountable spectra.

## 1. Introduction

A *theory* is a set of sentences - finite statements built from the function and relation symbols of a fixed language by the use of the Boolean connectives ("and", "not", etc.) and quantifiers ("there exists" and "for all"). The usual axioms for rings, groups and real closed fields are examples of theories. Associated to any theory is its class of *models*. A model of a theory is an algebraic structure that satisfies each of the sentences of the theory. For a theory $T$ and a cardinal $\kappa$, $I(T, \kappa)$ denotes the number of isomorphism classes of models of $T$ of size $\kappa$. The *uncountable spectrum* of a theory $T$ is the map $\kappa \mapsto I(T, \kappa)$, where $\kappa$ ranges over all uncountable cardinals. As examples, any theory of algebraically closed fields of fixed characteristic has $I(T, \kappa) = 1$ for all uncountable $\kappa$, while any completion of Peano Arithmetic has $I(T, \kappa) = 2^\kappa$.

With Theorem 6.1 we enumerate the possible uncountable spectra of complete theories in a countable language. Examples of theories possessing each of these spectra are given in [8].

Starting in 1970, Shelah placed the uncountable spectrum problem at the center stage of model theory. His goal was to develop a taxonomy of complete theories in a fixed countable language. Shelah's thesis was that the equivalence relation of 'having the same uncountable spectrum' induces a partition of the space of complete theories that is natural and useful for other applications. Over a span of twenty years he realized much of this research program. In addition to the results mentioned in the abstract, he showed that the uncountable spectrum $I(T, \kappa)$ is nondecreasing for all complete theories $T$ and that the divisions between spectra reflect structural properties. Shelah found a number of dividing lines among complete theories. The definitions of these dividing lines do not mention uncountable objects, but collectively they form an important distinction between the associated classes of uncountable models. On one hand, he showed that if a theory is on the 'nonstructure' side of at least one of these lines then models of the theory embed a certain amount of set theory; as a consequence their spectrum is maximal (i.e., $I(T, \kappa) = 2^\kappa$ for all uncountable $\kappa$). This is viewed as a negative feature, ruling out the possibility of a reasonable structure theorem for the class of models of the theory.

On the other hand, for models of theories that are on the 'structure' side of each of these lines, one can associate a system of combinatorial geometries. The isomorphism type of a model of such a theory is determined by local information (i.e., behavior of countable substructures) together with a system of numerical invariants (i.e., dimensions for the corresponding geometries). It follows that the uncountable spectrum of such a theory cannot be maximal. Thus, the uncountable spectrum of a complete theory in a countable language is nonmaximal if and only if every model of the theory is determined up to isomorphism by a well-founded, independent tree of countable substructures.

Our work is entirely contained in the stability-theoretic universe created by Shelah. We offer three new dividing lines on the space of complete theories in a (fixed) countable language (see Definitions 3.2 and 5.23) and show that these divisions, when combined with those offered by Shelah, are sufficient to characterize the uncountable spectra of all such theories. These new divisions partition the space of complete theories into Borel subsets (with respect to the natural topology on the space). The first two of these divisions measure how far the theory is from being totally transcendental, while the third division makes a much finer distinction between two spectra.

A still finer analysis in terms of geometric stability theory is possible. We mention for instance that any model of a complete theory whose uncountable spectrum is $\min\{2^{\aleph_\alpha}, \beth_{d-1}(|\alpha + \omega| + \beth_2)\}$ for some finite $d > 1$ interprets an infinite group. This connection turns out not to be needed for the present results, and will be presented elsewhere.

The main technical result of the paper is the proof of Theorem 3.3, which asserts that the embeddability of certain countable configurations of elements into some model of the theory gives strong lower bounds on the uncountable spectrum of the theory. The proof of this theorem uses techniques from descriptive set theory along with much of the technology developed to analyze superstable theories.

We remark that the computation of $I(T, \aleph_0)$ is still open. To wit, it remains unknown whether any countable, first-order theory $T$ has $I(T, \aleph_0) = \aleph_1$, even when the continuum hypothesis fails. Following our work, this instance is the only remaining open question regarding the possible values of $I(T, \kappa)$.

## 2. Background

Work on the spectrum problem is quite old. Morley's categoricity theorem [14], which asserts that if $I(T, \kappa) = 1$ for some uncountable cardinal $\kappa$, then $I(T, \kappa) = 1$ for all uncountable $\kappa$ is perhaps the most familiar computation of a spectrum. However, some work on the spectrum problem predates this. If $T$ has an infinite model, then for every infinite cardinal $\kappa$ the inequality $1 \leq I(T, \kappa) \leq 2^\kappa$ follows easily from the Löwenheim-Skolem theorem. Improving on this, Ehrenfeucht [4] discovered the notion of what is now called an unstable theory and showed that $I(T, \kappa) \geq 2$ for certain uncountable $\kappa$ whenever $T$ is unstable.

One cannot overemphasize the impact that Shelah has had on the the uncountable spectrum problem. Much of the creation of the subfield of stability theory is singlehandedly due to him and was motivated by this problem. The following survey of his definitions and theorems establish the framework for this paper and indicate why it is sufficient for us to work in the very restrictive setting of classifiable theories of finite depth.

Call a complete theory $T$ with an infinite model *classifiable* if it is superstable, has prime models over pairs, and does not have the dimensional order property. The following two theorems of Shelah indicate that this notion is a very robust dividing line among the space of complete theories.

THEOREM 2.1. *If a countable theory $T$ is not classifiable then $I(T, \kappa) = 2^\kappa$ for all $\kappa > \aleph_0$.*

*Proof.* If $T$ is not superstable then the spectrum of $T$ is computed in VIII 3.4 of [18]; this is the only place where this spectrum is computed for all uncountable cardinals $\kappa$. Hodges [9] contains a very readable proof of this for regular cardinals. Shelah's computation of the spectrum of a superstable theory with the dimensional order property is given in X 2.5 of [18]. More detailed proofs are given in Section 3 of Chapter XVI of Baldwin [1] and Theorem 2.3 of Harrington-Makkai [5]. Under the assumptions that $T$ is countable, superstable, and does not have the dimensional order property, the property of prime models over pairs is equivalent to $T$ not having the omitting types order property (OTOP). That the omitting types order property implies that $T$ has maximal spectrum was proved by Shelah in Chapter 12, Section 4 of [18]. Another exposition of this fact is given in [6]. □

In order to state the structural consequences of classifiability we state three definitions.

Definition 2.2.
1. $M$ is an *na-substructure* of $N$, $M \subseteq_{na} N$, if $M \subseteq N$ and for every formula $\varphi(x, y)$, tuple $a$ from $M$ and finite subset $F$ of $M$, if $N$ contains a solution to $\varphi(x, a)$ not in $M$ then $M$ contains a solution to $\varphi(x, a)$ that is not algebraic over $F$.

2. An $\omega$-*tree* $(I, <)$ is a partial order that is order-isomorphic to a nonempty, downward closed subtree of $^{<\omega}J$ for some index set $J$, ordered by initial segment. The root of $I$ is denoted by $\langle \rangle$ and for $\eta \neq \langle \rangle$, $\eta^-$ denotes the (unique) predecessor of $\eta$ in the tree.

3. An *independent tree of models* of $T$ is a collection $\{M_\eta : \eta \in I\}$ of elementary submodels of a fixed model of $T$ indexed by an $\omega$-tree $I$ that is independent with respect to the order on $I$.

4. A *normal tree of models* of $T$ is a collection $\{M_\eta : \eta \in I\}$ of models of $T$ indexed by an $\omega$-tree $I$ satisfying:

    - $\eta \lessdot \nu \lessdot \tau$ implies $M_\tau/M_\nu \perp M_\eta$;
    - for all $\eta \in I$, $\{M_\nu : \eta = \nu^-\}$ is independent over $M_\eta$.

5. A *tree decomposition* of $N$ is a normal tree of models $\{M_\eta : \eta \in I\}$ with the properties that, for every $\eta \in I$, $M_\eta$ is countable, $M_\eta \subseteq_{na} N$ and $\eta \lessdot \nu$ implies $M_\eta \subseteq_{na} M_\nu$ and $wt(M_\nu/M_\eta) = 1$.

THEOREM 2.3.  1. *Any normal tree of models is an independent tree of models.*

2. *If $T$ is classifiable then there is a prime model over any independent tree of models of $T$.*

3. *Every model $N$ of a classifiable theory is prime and minimal over any maximal tree decomposition contained in $N$.*

*Proofs.* The proof of (1) is an exercise in tree manipulations and orthogonality. Details can be found in Chapter 12 of [18] or Section 3 of Harrington-Makkai [5].

The proof of (2) only relies on $T$ having prime models over pairs. Its proof can be found in Chapter 12 of [18] or in [6].

The proof of (3) is more substantial. In [18] Shelah proves that any model of a classifiable theory has a number of tree decompositions of various sorts. However, the fact that a model of a classifiable theory admits a decomposition using countable, $na$-submodels is the content of Theorem C of Shelah-Buechler [19]. □

The two parts of Theorem 2.3 provide us with a method of producing upper bounds on spectra. Namely, $I(T,\kappa)$ is bounded above by the number of labelled $\omega$-trees of size $\kappa$. Since the components of the tree decompositions are countable, we may assume that the set of labels has size at most $2^{\aleph_0}$. In Section 5 we obtain better upper bounds in a number of cases by adding more structural information, which decreases the set of labels. However, at this point, there are still too many $\omega$-trees of size $\kappa$ to obtain a reasonable upper bound on $I(T,\kappa)$. A further reduction is available by employing the following additional definitions and theorems of Shelah.

*Definition* 2.4. An $\omega$-tree $(I, \triangleleft)$ is *well-founded* if it does not have an infinite branch. The *depth of a node* $\eta$ of a well-founded tree $I$ is defined inductively by
$$dp_I(\eta) = \sup\{dp_I(\nu) + 1 : \eta \triangleleft \nu\}$$
and the *depth of* $I$, $dp(I)$ is equal to $dp_I(\langle\rangle)$. A theory $T$ is *deep* if some model of $T$ has tree decomposition indexed by a non-well-founded tree. A (classifiable) theory $T$ is *shallow* if it is not deep. The *depth* $dp(T)$ of a shallow theory $T$ is the supremum of the depths of decomposition trees of models of $T$.

We remark that this definition of the depth of a theory differs slightly from the one given in [18]. The following theorem of Shelah is a consequence of Theorems X 5.1, X 4.7, and X 6.1 of [18]. Other proofs appear in Harrington-Makkai [5] and Baldwin [1].

THEOREM 2.5.   1. *If $T$ is classifiable and deep then $I(T, \kappa) = 2^\kappa$ for all $\kappa > \aleph_0$.*

2. *If $T$ is shallow then $dp(T) < \omega_1$ and, if $\omega \leq dp(T) < \omega_1$, then*
$$I(T, \aleph_\alpha) = \min\{2^{\aleph_\alpha}, \beth_{dp(T)+1}(|\alpha + \omega|)\}.$$

As a consequence of these results, we are justified in making the following assumption:

> *All theories in the rest of this paper are countable, classifiable and of finite depth d.*

For such theories, one obtains the naive upper bound of
$$I(T, \aleph_\alpha) \leq \beth_{d-1}(|\alpha + \omega|^{2^{\aleph_0}})$$
by induction on $d$, simply by counting the number of labelled trees in the manner described above.

In general, obtaining lower bounds on spectra is a challenging enterprise. The difficulty is due to the fact that the tree decomposition of a model given in Theorem 2.3 is typically not canonical. The method of *quasi-isomorphisms* introduced by Shelah and streamlined by Harrington-Makkai and Baldwin-Harrington (see e.g., [5, §3], [3, §3], or [1, Theorem XVII.4.7]) can be employed to show that if two models have 'sufficiently disparate' decomposition trees, then one can conclude that the models are nonisomorphic. From this, one obtains a (rather weak) general lower bound on the spectrum of a classifiable theory of finite depth $d > 1$, namely
$$I(T, \aleph_\alpha) \geq \min\{2^{\aleph_\alpha}, \beth_{d-2}(|\alpha + \omega|^{|\alpha+1|})\}.$$

A proof of this lower bound is given in Theorem 5.10(a) of [16].

Accompanying Shelah's 'top-down' analysis of the spectrum problem is work of Lachlan, Saffe, and Baldwin, who computed the spectra of theories satisfying much more stringent constraints.

In [11] and [12], Lachlan classifies the spectra of all $\omega$-categorical, $\omega$-stable theories. We use this classification verbatim at the end of Section 5. With [16], Saffe computes the uncountable spectra of all $\omega$-stable theories. A more detailed account of the analysis of this case is given by Baldwin in [1]. Aside from a few specific facts, we do not make use of this analysis here.

The history of this paper is modestly complicated. Shelah knew the value of $I(T, \aleph_\alpha)$ for large values of $\alpha$ (reported in Chapter XIII of [18]) modulo a certain continuum hypothesis-like question known as the SND (special number of dimensions) problem. In 1990, Hrushovski solved the SND problem; he also announced a calculation of the uncountable spectra. This calculation included a gap related to the behavior of nontrivial types but nonetheless a framework for the complete computation was introduced. The project lay fallow for several years before being taken up by the current authors; their initial work was reported in [8]. Hart and Laskowski recast the superstructure of the argument in a way to avoid the earlier gap and the work was completed while the three authors were in residence at MSRI.

We assume that the reader is familiar with stability theory. All of the necessary background can be found in the union of the texts by Baldwin [1] and Pillay [15]. Our notation is consistent with these texts. In Sections 3–6 we assume that we have a fixed, classifiable theory in a countable language of some finite depth. (The sufficiency of this assumption is explained in Section 2.) We work in $T^{eq}$. In particular, every type over an algebraically closed set is stationary. As well, throughout the paper we denote finite tuples of elements by singletons. For a stationary type $p(x)$ and a formula $\varphi(x, y)$, the notation $d_p x \varphi(x, y)$ denotes the $\varphi$-definition of $p$. As usual in the study of stable theories, we assume that we are working within the context of a large, saturated model $\mathfrak{C}$ of $T$. All sets are assumed to be subsets of $\mathfrak{C}$, and all models are assumed to be elementary submodels of $\mathfrak{C}$. In particular, the notation $M \subseteq N$ implies $M \preceq N$.

### 3. More dividing lines

In this section we provide a local analysis of a classifiable theory of finite depth $d$. In particular, we ask how far away the theory is from being totally transcendental. Towards this end, we make the following definitions.

*Definition* 3.1. 1. For any $n \leq dp(T)$, a *chain of length n*, $\overline{\mathcal{M}}$, is a sequence $M_0 \subseteq \ldots \subseteq M_{n-1}$ of $n$ countable models of $T$, where $M_{i+1}/M_i$ has weight 1 and $M_{i+1}/M_i \perp M_{i-1}$ for $i > 0$.

2. A chain $\overline{\mathcal{M}}$ is an *na-chain* if, in addition, each $M_i \subseteq_{na} \mathfrak{C}$.

3. For $\overline{\mathcal{M}}$ a chain of length $n$, the set of *relevant regular types* is the set
$$R(\overline{\mathcal{M}}) = \{p \in S(M_{n-1}) : p \text{ is regular and } p \perp M_{n-2}\}.$$
When $n = 1$, $R(\overline{\mathcal{M}})$ is simply the set of regular types over $M_0$.

4. A type $p \in R(\overline{\mathcal{M}})$ is *totally transcendental* (t.t.) *over* $\overline{\mathcal{M}}$ if there is a strongly regular $q \in R(\overline{\mathcal{M}})$, $q \not\perp p$ with a prime model $\overline{\mathcal{M}}(q)$ over $M_{n-1}$ and any realization of $q$.

It is clear that the notion of a relevant type being t.t. depends only on its nonorthogonality class. Our new dividing lines concern the presence or absence of a relevant type that fails to be t.t. and whether there is a trivial 'bad' type.

*Definition* 3.2.  A theory $T$ is *locally t.t. over* $\overline{\mathcal{M}}$ if every type in $R(\overline{\mathcal{M}})$ is t.t. over $\overline{\mathcal{M}}$. We say $T$ admits a *trivial failure* (of being t.t.) *over* $\overline{\mathcal{M}}$ if some trivial type $p \in R(\overline{\mathcal{M}})$ is not t.t. over $\overline{\mathcal{M}}$.

Our notation is consistent with standard usage, as any totally transcendental theory is locally t.t. over any chain. The heart of the paper will be devoted to showing that these conditions provide dividing lines for the spectra. In particular, the proof of the lower bounds offered below follows immediately from 3.17, 3.21, and 5.10.

THEOREM 3.3.   1. *If $T$ is not locally* t.t. *over some chain of length $n$ then*
$$I(T, \aleph_\alpha) \geq \begin{cases} \min\{2^{\aleph_\alpha}, \beth_2\} & \text{if } n = 1 \\ \min\{2^{\aleph_\alpha}, \beth_{n-2}(|\alpha + \omega|^{\beth_2})\} & \text{if } n > 1 \end{cases}$$
*for all ordinals $\alpha > 0$.*

2. *If $T$ admits a trivial failure over some chain of length $n$ then*
$$I(T, \aleph_\alpha) \geq \min\{2^{\aleph_\alpha}, \beth_{n-1}(|\alpha + \omega|^{\beth_2})\}$$
*for all ordinals $\alpha > 0$.*

Complementing this theorem, in Subsection 3.4 we show that if $T$ is locally t.t. over $\overline{\mathcal{M}}$, then there is a strong structure theorem for the class of weight-one models over $M_{n-1}$ that are orthogonal to $M_{n-2}$, especially when the chain $\overline{\mathcal{M}}$ has length equal to the depth of the theory.

The proof of Theorem 3.3 is broken into a number of steps. For the most part, the two parts of the theorem are proved in parallel. In Subsection 3.1 we define the crucial notions of diverse and diffuse families of leaves. Then, in the next two subsections we analyze the two ways in which a theory could fail to be locally t.t. over a chain. For each of these, we will show that the failure of being locally t.t. over some chain of length $n$ implies the existence of a diverse family of leaves of size continuum over some $na$-chain of length $n$. In addition, if there is a trivial failure of being t.t., then the family of leaves mentioned above will actually be diffuse. Then, in Section 4, we establish some machinery to build many nonisomorphic models from the existence of a diverse or diffuse family of leaves. Much of this is bookkeeping, but there are two important ideas developed there. Foremost is the Unique Decomposition Theorem (Theorem 4.1), which enables us to preserve nonisomorphism as we step down a decomposition tree. The other idea which is used is the fact that decomposition trees typically have many automorphisms. This fact implies that models that are built using such trees as skeletons have desirable homogeneity properties (see Definition 4.3). Finally, we complete the proof of Theorem 3.3 in Section 5.

3.1. *Diverse and diffuse families of leaves.* We begin by introducing some convenient notation for prime models over independent trees of models. First of all, suppose that $N_1$ and $N_2$ are two submodels of our fixed saturated model which are independent over a common submodel $N_0$. By $N_1 \oplus_{N_0} N_2$ we will mean a prime model over $N_1 \cup N_2$; this exists because we are assuming our theory has prime models over pairs. For our purposes, the exact model that we fix will not matter because we are only interested in its isomorphism type. In abstract algebraic terms, we want to think of this as an "internal" direct sum.

Now suppose that $M_0$ is any model and $M_0 \subseteq M_i$ for $i = 1, 2$, not necessarily independent over $M_0$. By $M_1 \oplus_{M_0} M_2$ we will mean the "external" direct sum i.e., we choose $M_i'$ isomorphic to $M_i$ over $M_0$ and such that $M_1'$ is independent over $M_2'$ over $M_0$ and form $M_1' \oplus_{M_0} M_2'$ in the internal sense defined above. We similarly define $\bigoplus_M \mathcal{F}$ for a family of models, each of which contains a fixed model $M$.

Suppose that $\mathcal{M} = \langle M_\zeta : \zeta \in I \rangle$ is an independent family of models with respect to a tree ordering $\langle I, \lessdot \rangle$ such that if $\eta \lessdot \zeta \in I$ then $M_\eta \subseteq M_\zeta$. We say that a family of elementary maps $\langle f_\zeta : \zeta \in I \rangle$ is *compatible* with $\mathcal{M}$ if whenever $\zeta \in I$ we have $\mathrm{dom}(f_\zeta) = M_\zeta$ and if $\eta \lessdot \zeta \in I$ then $f_\zeta|_{M_\eta} = f_\eta$.

Definition 3.4. If $\overline{\mathcal{M}}$ is an $na$-chain of length $n$, then the *set of leaves of $\overline{\mathcal{M}}$*, Leaves($\overline{\mathcal{M}}$), is a set containing one representative up to isomorphism over $\overline{\mathcal{M}}$ of all chains $\overline{\mathcal{N}}$ of length $(n+1)$ extending $\overline{\mathcal{M}}$. If $\overline{\mathcal{M}}$ is an $na$-chain of length $n$ and $Y \subseteq$ Leaves($\overline{\mathcal{M}}$), then an $(\overline{\mathcal{M}}, Y)$-*tree* is an independent tree of models $\mathcal{M} = \langle M_\zeta : \zeta \in I \rangle$ where $I$ has height at most $n+1$ together with a distinguished copy of $\overline{\mathcal{M}}$ and a family of elementary maps $\langle f_\zeta : \zeta \in I \rangle$ compatible with $\mathcal{M}$ such that $f_\zeta$ maps $M_\zeta$ onto $N_{lg(\zeta)}$ for some $\overline{\mathcal{N}} \in Y$.

(In particular, note that if $lg(\zeta) \leq n$ then $f_\zeta$ maps onto $M_{lg(\zeta)}$.) An $(\overline{\mathcal{M}}, Y)$-*model* is a model which is prime over an $(\overline{\mathcal{M}}, Y)$-tree; the copy of $\overline{\mathcal{M}}$ in this tree will be distinguished as a chain of submodels of this $(\overline{\mathcal{M}}, Y)$-model.

We make an important convention: Suppose $N_1$ and $N_2$ are two $(\overline{\mathcal{M}}, Y)$-models and we wish to form $N_1 \oplus_{M_k} N_2$ for some $k < n$. We declare that this sum is an "external" sum as discussed above. If we wish to view this new model as an $(\overline{\mathcal{M}}, Y)$-model, we need to specify which copy of $\overline{\mathcal{M}}$ will be distinguished in the sum. Our convention is that the distinguished copy in the sum will be the distinguished copy from the left-most summand.

As notation, if $Z$ is a set of $\overline{\mathcal{M}}$-leaves then

$$N^*(Z) = \bigoplus_{M_{n-1}} \{N : \overline{N} \in Z \text{ where } \overline{N}(n) = N\}.$$

We next isolate the two crucial properties of a set $Y$ of leaves. In Section 5 we will show the effects on lower bound estimates for spectra given that there are large families of leaves with these properties. Lemma 4.2 of Section 4 will show that a diffuse family is diverse.

*Definition* 3.5. Let $\overline{\mathcal{M}}$ be an $na$-chain of length $n$.

1. A set $Y \subseteq \text{Leaves}(\overline{\mathcal{M}})$ is *diffuse* if

$$N \oplus_{M_{n-1}} V \not\cong_V N' \oplus_{M_{n-1}} V$$

for all distinct $N, N' \in Y$ and any $(\overline{\mathcal{M}}, Y)$-model $V$.

2. A set $Y \subseteq \text{Leaves}(\overline{\mathcal{M}})$ is *diverse* if

$$N^*(Z_1) \oplus_{M_{n-2}} V \not\cong N^*(Z_2) \oplus_{M_{n-2}} V$$

over $V$ and $M_{n-1}$ for all distinct $Z_1, Z_2 \subseteq Y$ and any $(\overline{\mathcal{M}}, Y)$-model $V$.

If $n = 1$ we omit the model $V$ and the condition becomes: for distinct $Z_1, Z_2 \subseteq Y$,

$$N^*(Z_1) \not\cong_{M_0} N^*(Z_2)$$

The following lemma provides us with an easy way of producing a diffuse family of leaves.

LEMMA 3.6. *Suppose that $\{p_i : i \in I\}$ is a set of pairwise orthogonal types in $R(\overline{\mathcal{M}})$, and that $\{N_i : i \in I\} \subseteq \text{Leaves}(\overline{\mathcal{M}})$ is a set of models such that $N_i$ is dominated by a realization of $p_i$ over $M_{n-1}$. Then $\{N_i : i \in I\}$ is diffuse.*

*Proof.* In fact, for any model $V$ containing $M_{n-1}$, if $N_i \oplus_{M_{n-1}} V \cong_V N_j \oplus_{M_{n-1}} V$ then $p_i$ is not orthogonal to $p_j$ so $i = j$. □

We conclude this subsection by introducing the notion of a *special* type, showing that they are present in every nonorthogonality class of $R(\overline{\mathcal{M}})$, and proving two technical lemmas that will be used in subsections 3.2 and 3.3 to obtain diverse or diffuse families of leaves.

As notation for a regular strong type $p$, let $[p]$ denote the collection of strong types (over any base set) nonorthogonal to $p$. We let $R^\infty([p]) = \min\{R^\infty(q) : q \in [p]\}$. It is easy to see that $R^\infty([p])$ is the smallest ordinal $\alpha$ such that $p$ is nonorthogonal to some formula $\theta$ of $R^\infty$-rank $\alpha$, which is also the smallest ordinal $\beta$ such that $p$ is foreign to some formula $\psi$ of $R^\infty$-rank $\beta$. The following lemma is general and holds for any superstable theory.

LEMMA 3.7. *Let $M^- \subseteq M$ be models of a superstable theory and suppose that a regular type $p \in S(M)$ is orthogonal to $M^-$ but is nonorthogonal to some $\theta(x, b)$, where $R^\infty(\theta(x, b)) = R^\infty([p])$. Then $\sigma(p)$ is foreign to $\theta(x, b)$ for any automorphism $\sigma \in \mathrm{Aut}_{M^-}(\mathfrak{C})$ satisfying $\sigma(M) \underset{M^-}{\downarrow} b$.*

*Proof.* Suppose that $\theta(c, b)$ holds and let $X$ be any set. We claim that $\mathrm{tp}(c/Xb) \perp \sigma(p)$. To see this, first note that $\sigma(p) \perp M^- b$, since we assumed $\sigma(p) \perp M^-$ and $\sigma(M) \underset{M^-}{\downarrow} b$. There are now two cases. If $R^\infty(c/Xb) = R^\infty([p])$, then $\mathrm{tp}(c/Xb)$ does not fork over $b$, so $\sigma(p) \perp \mathrm{tp}(c/Xb)$ by the note above. On the other hand, if

$$R^\infty(c/Xb) < R^\infty([p]) = R^\infty([\sigma(p)]),$$

then $\mathrm{tp}(c/Xb) \perp \sigma(p)$ as well. □

We now turn our attention back to a particular chain $\overline{\mathcal{M}}$ of length $n$.

*Definition* 3.8.   1. If $p \in R(\overline{\mathcal{M}})$ then $q$ is a *tree conjugate* of $p$ if for some $k < n - 1$ there is an automorphism $\sigma$ fixing $M_k$ pointwise such that $\sigma(p) = q$ and $\sigma(M) \underset{M_k}{\downarrow} M_{n-1}$. (If $n = 1$ then $p$ does not have any tree conjugates.)

2. A type $p \in R(\overline{\mathcal{M}})$ is *special via* $\varphi(x, e)$ if $\varphi(x, e)$ is $p$-simple, $\varphi(x, e) \in p$, and the tree conjugates of $p$ are foreign to $\varphi(x, e)$. A type $p \in R(\overline{\mathcal{M}})$ is *special* if it is special via some formula.

As promised, we show that special types exist in every nonorthogonality class of $R(\overline{\mathcal{M}})$. The proof of the following lemma is an adaptation of the proof of Lemma 8.2.19 of [15], which in turn is adapted from arguments in [18].

LEMMA 3.9. *Let $a$ be any realization of a type $p \in R(\overline{\mathcal{M}})$, where $\overline{\mathcal{M}}$ is a chain of length $n$. There is an $a' \in \mathrm{acl}(M_{n-1}a) \smallsetminus M_{n-1}$ such that $p' = \mathrm{tp}(a'/M_{n-1}) \in R(\overline{\mathcal{M}})$ is special.*

*Proof.* Choose a formula $\theta(x,b)$ nonorthogonal to $p$ with $R^\infty(\theta(x,b)) = R^\infty([p])$ and choose a (regular) type $q$ nonorthogonal to $p$ containing $\theta(x,b)$. Choose a set $A \supseteq M_{n-1}$ and a nonforking extension $r \in S(A)$ of $q$ such that $a \underset{M_{n-1}}{\downarrow} A$, $r$ is stationary and there is a realization $c$ of $r$ with $a \underset{A}{\downarrow} c$. Now choose $a' \in Cb(\mathrm{stp}(bc/Aa)) \smallsetminus \mathrm{acl}(A)$. Since $\mathrm{tp}(a/A)$ does not fork over $M_{n-1}$, $a' \in \mathrm{acl}(M_{n-1}a) \smallsetminus M_{n-1}$, hence $p' = \mathrm{tp}(a'/M_{n-1})$ is regular and nonorthogonal to $p$. It remains to find an $L(M_{n-1})$-formula witnessing that $\mathrm{tp}(a'/M_{n-1})$ is special. Choose a Morley sequence $I = \langle b_n c_n : n \in \omega \rangle$ in $\mathrm{stp}(bc/Aa)$ with $b_0 c_0 = bc$. Since $a' \in \mathrm{dcl}(\overline{b}\overline{c})$ for some initial segment of $I$, $a' = f(\overline{b}, \overline{c})$ for some $\emptyset$-definable function $f$. Note that $a \underset{M_{n-1}}{\downarrow} A\overline{b}$. Choose a finite $A_0 \subseteq A$ on which everything is based, let $w = \mathrm{tp}(A_0 \overline{b}/M_{n-1})$ and let $\varphi(x,e)$ be the $L(M_{n-1})$-formula

$$\varphi(x,e) := d_w \overline{y} \exists \overline{z} \left( \bigwedge_i \theta(z_i, y_i) \wedge x = f(\overline{y}, \overline{z}) \right).$$

For notation, assume that $w$ is based on $e$. Clearly $\varphi(a', e)$ holds and it follows easily that $\varphi(x,e)$ is $p'$-simple. Hence if $n = 1$ we are done. So suppose that $n > 1$ and $\varphi(a^*, e)$ holds. Fix $k < n-1$ and an automorphism $\sigma$ over $M_k$ such that $\sigma(M_{n-1}) \underset{M_k}{\downarrow} M_{n-1}$. Choose $A^* \overline{b}^*$ realizing $w|M_{n-1}\sigma(M_{n-1})$ and choose $\overline{c}^*$ such that $a^* = f(\overline{b}^*, \overline{c}^*)$ and $\theta(c_i^*, b_i^*)$ holds for each $i$. Since $p'$ is not orthogonal to $p$, $\theta(x, b_i^*)$ is nonorthogonal to $p'$ and has least $R^\infty$-rank among all formulas nonorthogonal to $p'$. Thus, as $b_i^* \underset{M_k}{\downarrow} \sigma(M_{n-1})$, it follows from Lemma 3.7 that $\sigma(p')$ is foreign to $\theta(x, b_i^*)$ for all $i$. Hence $\mathrm{tp}(a^*/e)$ is hereditarily orthogonal to $\sigma(p')$. That is, the tree conjugates of $p'$ are foreign to $\varphi(x,e)$. □

The following lemma is simply a restatement of Lemma 3.9 together with an application of the Open Mapping Theorem.

LEMMA 3.10. *Suppose that $\overline{\mathcal{M}}$ is a chain of length $n > 1$ and $p \in R(\overline{\mathcal{M}})$ is special via $\varphi$. Further, suppose that a model $N/M_{n-1} \perp M_{n-2}$, and $U$ is dominated over $W$ by an independent set of conjugates of $p$. Then any realization of $\varphi$ in $N \oplus_{M_{n-2}} U$ is contained in $N \oplus_{M_{n-2}} W$.*

*Proof.* Suppose that $a \in N \oplus_{M_{n-2}} U$ realizes $\varphi$. Then $\mathrm{tp}(a/NU)$ is isolated. As the tree conjugates of $p$ are foreign to $\varphi$,

$$a \underset{NW}{\downarrow} d$$

for any $d \in U$ realizing a conjugate of $p$. Thus $a \underset{NW}{\downarrow} U$, since $U$ is dominated over $W$ by an independent set of realizations of conjugates of $p$. Hence $tp(a/NW)$ is isolated, which implies that $a \in N \oplus_{M_{n-2}} W$. □

If, in addition, our special type is trivial then we can say more. The lemma that follows is one of the main reasons why we are able to build a diffuse family instead of a diverse family when the 'offending' type is trivial.

LEMMA 3.11. *Suppose that $\overline{\mathcal{M}}$ is a chain of length $n$, $N \in \text{Leaves}(\overline{\mathcal{M}})$, $p \in R(\overline{\mathcal{M}})$ is any trivial, special type via $\varphi$ such that some realization of $p$ dominates $N$ over $M_{n-1}$. Suppose further that $M_{n-1} \subseteq W \subseteq U$, where $U$ is dominated by $W$-independent realizations of nonforking extensions of tree conjugates of $p$ over $W$ and regular types not orthogonal to $p$. Then if an element $a$ satisfies $\varphi$, $a \in N \oplus_{M_{n-1}} U$, $tp(a/M_{n-1})$ is regular and $a \underset{M_{n-1}}{\downarrow} U$ then $tp(a/NW)$ is isolated.*

*Proof.* Note first that the assumptions imply that $tp(a/M_{n-1})$ is not orthogonal to $p$. Since $p$ is trivial and $a \underset{M_{n-1}}{\downarrow} U$, it follows that $a$ forks with $N$ over $M_{n-1}$. Since $\varphi(a)$ holds, it follows that $tp(a/N)$ (and any extension of this type) is orthogonal to $p$ and all tree conjugates of $p$. It follows that $a \underset{NW}{\downarrow} U$. Since $tp(a/NU)$ is isolated, it follows that $tp(a/NW)$ is isolated. □

3.2. *The existence of prime models.* Throughout this section we assume that there is some chain $\overline{\mathcal{M}}$ of length $n$ together with a type $r \in R(\overline{\mathcal{M}})$ for which there is no prime model over $M_{n-1}c$, where $c$ is a realization of $r$. By choosing an extension $\overline{\mathcal{M}}'$ of $\overline{\mathcal{M}}$ which is minimal in a certain sense, we will construct a highly disparate family $Y = \{N_\eta : \eta \in {}^\omega 2\}$ of $\text{Leaves}(\overline{\mathcal{M}}')$ and a family of types $\{s_\eta(x, z)\}$ over $M'_{n-1}$ that will witness this disparity. Then, following the construction of the family in Proposition 3.14, we argue in Corollary 3.17 that if the original type was special, then this set of leaves is diverse. Further, if in addition the type $r$ is trivial, we show that this family is actually diffuse. We begin by specifying what we mean by a free extension of a chain.

Definition 3.12. An $na$-chain $\overline{\mathcal{M}}'$ *freely extends* the chain $\overline{\mathcal{M}}$ if both chains have the same length (say $n$), $M_0 \subseteq M'_0$, $M'_{i-1} \cup M_i \subseteq M'_i$, and $M'_{i-1} \underset{M_{i-1}}{\downarrow} M_i$ for all $0 < i < n$.

It is readily checked that if $r \in R(\overline{\mathcal{M}})$ is special, then for any free extension $\overline{\mathcal{M}}'$ of $\overline{\mathcal{M}}$, the nonforking extension of $r$ to $R(\overline{\mathcal{M}}')$ will be special as well. The bulk of this subsection is devoted to the proof of Proposition 3.14. In order to state the proposition precisely, we require some notation.

Suppose that $Y$ is a family of Leaves($\overline{\mathcal{M}}$) that is indexed by $^\omega 2$. Fix an $(\overline{\mathcal{M}}, Y)$-model $V$ and a sequence $\eta \in {}^\omega 2$. We can decompose $V$ into two pieces,

$$V = V_\eta \bigoplus_{W_V} V_{\text{no } \eta}$$

where $W_V$ is the model prime over the tree truncated below level $n$,

$$V_\eta = \bigoplus_{W_V} \{N_i : N_i \text{ conjugate to } N_\eta\}$$

and

$$V_{\text{no } \eta} = \bigoplus_{W_V} \{N_i : N_i \text{ not conjugate to } N_\eta\}.$$

*Definition* 3.13. A formula $\theta$ is $\psi$-*definable* over $A$ if $\theta(\mathfrak{C}) \subseteq \text{dcl}(A \cup \psi(\mathfrak{C}))$.

The following special case will be used in the proof of the proposition below: If $M \subseteq A$, $\theta \in L(A)$, $\psi \in L(M)$ and for every $a$ realizing $\theta$ there is $b$ such that $b \underset{M}{\downarrow} Aa$ and $a \in \text{dcl}(Ab \cup \psi(\mathfrak{C}))$, then it follows easily from compactness and the finite satisfiability of nonforking over a model that $\theta$ is $\psi$-definable over $A$.

PROPOSITION 3.14. *Assume that $\overline{\mathcal{M}}$ is a chain of length $n$ and $r \in R(\overline{\mathcal{M}})$ is a type such that there is no prime model over $M_{n-1}c$ for $c$ a realization of $r$. Then there is a free extension $\overline{\mathcal{M}}'$ of $\overline{\mathcal{M}}$ and a family $Y = \{N_\eta : \eta \in {}^\omega 2\}$ of Leaves($\overline{\mathcal{M}}'$), along with a family $\{s_\eta(x, z) : \eta \in {}^\omega 2\}$ of types over $M'_{n-1}$ such that each $N_\eta$ realizes $s_\eta$, yet for any $(\overline{\mathcal{M}}', Y)$-model $V$, $V$ omits $s_\eta(x, c^*)$ for all $c^* \in V_{\text{no } \eta}$ where $c^*$ realizes $r|M'_{n-1}$.*

*Proof.* The lack of a prime model over $M_{n-1}c$ implies the lack of a prime model over $\text{acl}(M_{n-1}c)$. Look at all possible quadruples $(\overline{\mathcal{M}}', r', \theta, \psi)$ where $\overline{\mathcal{M}}'$ is a free extension of $\overline{\mathcal{M}}$, $r'$ is the nonforking extension of $r$ to $M'_{n-1}$, and (fixing a realization $c$ of $r'$ and letting $C = \text{acl}(M'_{n-1}c)$) $\theta(x)$ is a formula in $L(C)$ with no isolated extensions over $C$ and $\psi$ is a formula in $L(M'_{n-1})$ such that $\theta(x)$ is $\psi$-definable over $C$. Among all such quadruples, fix one for which $R^\infty(\psi)$ is minimal. To ease notation in what follows, we denote $\overline{\mathcal{M}}'$ by $\overline{\mathcal{M}}$ and $r'$ by $r$. As well, fix a realization $c$ of $r$ and let $C = \text{acl}(M_{n-1}c)$.

Let $L(C)^*$ denote the set of consistent $L(C)$-formulas that are $R^\infty$-minimal over $C$. (See Definition B.3 and the discussion following it for the utility of restricting to this class of formulas.) It is easily seen that every consistent $L(C)$-formula extends to an $R^\infty$-minimal formula over $C$ and that every consistent extension of such a formula remains $R^\infty$-minimal over $C$. In particular, we may assume that $\theta$ is $R^\infty$-minimal over $C$.

We will construct the families of leaves and types simultaneously by building successively better finite approximations. We adopt the notation of forcing (i.e., partial orders and filters meeting collections of dense sets). However, as we will only insist that our filters meet a countable collection of dense sets, a 'generic object' will already be present in the ground model and each of these constructions could just as easily be considered as a Henkin construction. Our set of forcing conditions $\mathcal{P}$ is the set of functions $p$ that satisfy the following conditions:

1. $\mathrm{dom}(p) = {}^{\leq m}2$ for some $m \in \omega$;

2. $p(\eta) \in L(C)^*$ and has its free variables among $\{x_k^\eta : k \in \omega\}$;

3. If $\nu, \eta \in \mathrm{dom}(p)$, $\nu \triangleleft \eta$, and $\varphi \in L(C)^*$, then $p(\nu) \vdash \varphi(x_0^\nu, \ldots, x_{k-1}^\nu)$ implies $p(\eta) \vdash \varphi(x_0^\eta, \ldots, x_{k-1}^\eta)$; and

4. $p(\langle\rangle) \vdash \theta(x_0^{\langle\rangle})$.

We let $m(p)$ be the $m$ such that $\mathrm{dom}(p) = {}^{\leq m}2$. If $p, q \in \mathcal{P}$ we put $p \leq q$ if and only if $m(p) \geq m(q)$ and $p(\eta) \vdash q(\eta)$ for all $\eta \in \mathrm{dom}(q)$. If $G \subseteq \mathcal{P}$ is any filter and $\eta \in {}^\omega 2$, let

$$p_\eta(G) = \{\theta \in L(C) : p(\nu) \vdash \theta(x_0^\nu, \ldots, x_{k-1}^\nu) \text{ for some } p \in G \text{ and some } \nu \triangleleft \eta\}.$$

We first list a set of basic density conditions we want our filter to meet:

1. For each $\varphi \in L(C)$ and each $k \in \omega$,

$$D_{\varphi,k} = \{p \in \mathcal{P} : m(p) \geq k \text{ and } \forall \eta \in {}^{m(p)}2(p(\eta) \vdash \varphi \text{ or } p(\eta) \vdash \neg\varphi)\};$$

2. For each $\psi(y, z) \in L(C)$ and each $k \in \omega$, let $D_{\psi,k}$ be

$$\{p \in \mathcal{P} : m(p) \geq k \text{ and } \forall \eta \in {}^{m(p)}2(p(\eta) \vdash \neg\exists z \psi(y, z) \vee \psi(y, x_j^\eta) \text{ for some } j)\}.$$

It is easy to see that every basic condition is dense in $\mathcal{P}$. As well, if $G$ is a filter meeting all of the conditions mentioned above (i.e., $G \cap D \neq \emptyset$ for each $D$) then for each $\eta \in {}^\omega 2$, $p_\eta(G)$ is a complete type (in an $\omega$-sequence of variables) over $C$ and if $\langle b_k : k \in \omega \rangle$ is a realization of $p_\eta(G)$ then $N_\eta = \{b_k : k \in \omega\}$ is a leaf of $\overline{\mathcal{M}}$ that contains $C$. (The fact that $wt(N_\eta/M_{n-1}) = 1$ and $N_\eta/M_{n-1} \perp M_{n-2}$ follows from our choice of $L(C)^*$ and Lemma B.4.) Additionally, $N_\eta \models \theta(b_0, c)$. In what follows, we will take $s_\eta(x, z)$ to be $\mathrm{tp}(b_0 c/M_{n-1})$.

Before stating the crucial density conditions that will ensure that the families of leaves and types satisfy the conclusions of the proposition, we pause to set notation. If $I$ is a tree in which every branch has length $n + 1$, let $I^+ = \{\zeta \in I : lg(\zeta) = n\}$.

*Definition* 3.15.   A *finite approximation* $\mathcal{F}$ of height $m$ consists of a finite independent tree $\langle B_\zeta : \zeta \in I \rangle$ of sets where every branch of $I$ has length $n + 1$, together with a distinguished copy of $\overline{\mathcal{M}}$, a family of elementary maps $\langle f_\zeta : \zeta \in I \rangle$ compatible with $\overline{\mathcal{M}}$, and a map $\pi : I^+ \to {}^m 2$. We require that $f_\zeta$ maps $M_{lg(\zeta)}$ onto $B_\zeta$ if $lg(\zeta) < n$ and $f_\zeta$ maps $C$ onto $B_\zeta$ if $lg(\zeta) = n$.

As notation, $B^* = \bigcup \{B_\zeta : \zeta \in I\}$, and $\mathrm{Var}(\mathcal{F}) = \{x_\zeta^n : n \in \omega, \zeta \in I^+\}$. We say that the finite approximation $\mathcal{F}'$ is a *natural extension* of $\mathcal{F}$ if the sets $I$, $\langle B_\zeta : \zeta \in I \rangle$, the choice of $\overline{\mathcal{M}}$, and the maps $\langle f_\zeta : \zeta \in I \rangle$ of $\mathcal{F}$ and $\mathcal{F}'$ are identical; the height of $\mathcal{F}'$ is at least the height of $\mathcal{F}$, and $\pi'(\zeta) \triangleleft \pi(\zeta)$ for all $\zeta \in I^+$.                                                                                                        □

Suppose $p \in \mathcal{P}$ and $\mathcal{F}$ is a finite approximation of height $m(p)$. Let

$$\mathcal{F}(p) = \bigwedge_{\zeta \in I^+} f_\zeta(p(\pi(\zeta))).$$

The formula $\mathcal{F}(p)$, which is over $B^*$ and whose free variables are among $\mathrm{Var}(\mathcal{F})$, should be thought of as an approximation to an $(\overline{\mathcal{M}}, Y)$-model $V$ in the statement of the proposition.

As notation, for a fixed finite approximation $\mathcal{F}$ of height $m$ and $\eta \in {}^m 2$, let $I_\eta^+ = \{\zeta \in I^+ : \pi(\zeta) = \eta\}$, and let $\mathrm{Var}_\eta(\mathcal{F}) = \{x_\zeta^n : n \in \omega, \zeta \in I_\eta^+\}$. Let $I_{\mathrm{no}\ \eta} = I \setminus I_\eta^+$, $B^*_{\mathrm{no}\ \eta} = \bigcup \{B_\zeta : \zeta \in I_{\mathrm{no}\ \eta}\}$ and $\mathrm{Var}_{\mathrm{no}\ \eta}(\mathcal{F}) = \mathrm{Var}(\mathcal{F}) \setminus \mathrm{Var}_\eta(\mathcal{F})$.

We now introduce the crucial set of density conditions.

*Density condition* 3.16.   Fix a finite approximation $\mathcal{F}$ of height $m$ and an $\eta \in {}^m 2$. Choose an $L(B^*_{\mathrm{no}\ \eta})$-formula $\chi(y, \overline{u})$ with $\overline{u} \subseteq \mathrm{Var}_{\mathrm{no}\ \eta}(\mathcal{F})$ and an $L(B^*)$-formula $\varphi(x, y, \overline{v})$ such that $\overline{u} \subseteq \overline{v} \subseteq \mathrm{Var}(\mathcal{F})$. Let $D = D(\mathcal{F}, \eta, \chi, \varphi)$ consist of all $p \in \mathcal{P}$ such that $m(p) \geq m$ and for some natural extension $\mathcal{F}'$ of $\mathcal{F}$ of height $m(p)$,

1. $\mathcal{F}'(p) \vdash \forall y \neg \chi(y, \overline{u}) \vee \exists y(\chi(y, \overline{u}) \wedge \neg \delta(y))$ for some $\delta \in r$ or

2. $\mathcal{F}'(p) \vdash \neg \exists x \exists y(\chi(y, \overline{u}) \wedge \varphi(x, y, \overline{v}))$ or

3. $\mathcal{F}'(p) \vdash \exists x \exists y(\chi(y, \overline{u}) \wedge \varphi(x, y, \overline{v}) \wedge \neg \delta(x, y))$ for some $L(M_{n-1})$-formula $\delta$ satisfying $p(\eta) \vdash \delta(x_0^\eta, c)$.

We check that if $G$ is a filter that meets the basic conditions and intersects each of the sets $D(\mathcal{F}, \eta, \chi, \varphi)$, then the sets of leaves $Y = \{N_\eta : \eta \in {}^\omega 2\}$ and types $\{s_\eta(x, z) : \eta \in {}^\omega 2\}$ satisfy the conclusion of the proposition.

Toward this end, fix an $(\overline{\mathcal{M}}, Y)$-tree, an $\eta \in {}^\omega 2$ and suppose $V$ (respectively $V_{\mathrm{no}\ \eta}$) is prime over this tree (respectively over the leaves not conjugate to $N_\eta$). Further, suppose by way of contradiction that $c^* \in V_{\mathrm{no}\ \eta}$ realizes $r$ and $b \in V$ realizes $s_\eta(x, c^*)$. Now in fact $c^*$ and $b$ are isolated over a finite part of the given $(\overline{\mathcal{M}}, Y)$-tree; $c^*$ is isolated over this tree by a formula $\chi(y)$ and $b$ is isolated over $c^*$ and the tree by a formula $\varphi(x, c^*)$. We suppress the parameters from the tree to ease notation.

Choose a number $m$ such that if conjugates of $N_\nu$ and $N_\mu$ appear in the finite tree needed to isolate $c^*$ and $b$ then if $\nu|_m = \mu|_m$ then $\nu = \mu$. Now this finite tree together with $\eta$, $\chi$ and $\varphi$ lead naturally to a finite approximation $\mathcal{F}$ of height $m$ and a density condition $D = D(\mathcal{F}, \eta|_m, \chi, \varphi)$. We claim that $D$ cannot meet the filter $G$. Choose any $p \in G$. By replacing $\mathcal{F}$ by a natural extension, we may assume $m(p) = m$. Now from the conditions for $D$, clearly the first condition must have failed since $\chi(y)$ is consistent and implies $r$. The second condition must also have failed because $\varphi(x, y) \wedge \chi(y)$ is consistent. But, if the third condition held, then $\varphi(x, c^*)$ would not even isolate $tp(b/M_{n-1}c^*)$. Hence $D \cap G = \emptyset$.

Thus, to complete the proof, it remains to show that each of the sets $D = D(\mathcal{F}, \eta, \chi, \varphi)$ is dense in $\mathcal{P}$. So fix $D$ and choose $p \in \mathcal{P}$. Without loss, $\mathcal{F}$ has height $m = m(p)$.

By an $\mathcal{F}$-*potential extension of* $p$ we mean a sequence of types $\bar{q} = \langle q_\nu : \bar{\nu} \in {}^m 2 \rangle$ such that each $q_\nu \in S(C)$ extends $p(\nu)$, together with a realization $\bar{a}$ of

$$\mathcal{F}(\bar{q}) = \bigwedge_{\zeta \in I^+} f_\zeta(q_{\pi(\zeta)}).$$

As each $q_\nu$ is $c$-isolated, the set of elements of $\bar{a}$ are independent over $B^*$, hence $\mathcal{F}(\bar{q})$ is a complete type over $B^*$. We concentrate on three cases which correspond to the three conditions in Density Condition 3.16.

*Case* one. Does there exist an $\mathcal{F}$-potential extension $(\bar{q}, \bar{a})$ of $p$ such that $\neg \exists y \chi(y, \bar{a}) \vee \exists y (\chi(y, \bar{a}) \wedge \neg \delta(y))$ holds for some $\delta \in r$? If so, then by Lemma B.4(3) we can find a sequence of $L(C)$-formulas $\langle \alpha_\nu : \nu \in {}^m 2 \rangle$ such that each $\alpha_\nu \in q_\nu$ extends $p(\nu)$ so that if we define $p' \leq p$ by

$$p'(\nu) = \begin{cases} p(\nu) & \text{if } \nu \in {}^{<m} 2; \\ \alpha_\nu & \text{if } \nu \in {}^m 2, \end{cases}$$

then $\mathcal{F}(p') \vdash \neg \exists y \chi(y, \overline{u}) \vee \exists y (\chi(y, \bar{u}) \wedge \neg \delta(y))$ for some $\delta \in r$.

*Case* two. Does there exist an $\mathcal{F}$-potential extension $(\bar{q}, \bar{a})$ of $p$ such that $\neg \exists x \exists y (\chi(y, \bar{a}) \wedge \varphi(x, y, \bar{a}))$ holds? If so, as in the first case we can use Lemma B.4(3) to define $p' \leq p$ such that $\mathcal{F}(p') \vdash \neg \exists x \exists y (\chi(y, \bar{u}) \wedge \varphi(x, y, \bar{v}))$.

*Case* three. Does there exist an $\mathcal{F}$-potential extension $(\bar{q}, \bar{a})$ of $p$ such that $\exists x \exists y (\chi(y, \bar{a}) \wedge \varphi(x, y, \bar{a}) \wedge \neg \delta(x, y))$ holds for some $\delta(x, y) \in L(M_{n-1})$ such that $\delta(x_0, c) \in q_\eta$? If so, then using Lemma B.4(3) we can define $p' \leq p$ so that $\mathcal{F}(p') \vdash \exists x \exists y (\chi(y, \bar{u}) \wedge \varphi(x, y, \bar{v}) \wedge \neg \delta(x, y))$.

We complete the proof that $D$ is dense by showing that these three cases are exhaustive. So, assume by way of contradiction that for any $\mathcal{F}$-potential extension $(\overline{q}, \overline{a})$ of $p$, $\chi(y, \bar{a})$ is consistent and implies $r(y)$, and $\chi(y, \bar{a}) \wedge \varphi(x, y, \bar{a})$ is consistent and implies $q_\eta(x, y) \in L(M_{n-1})$. We will show that this is an impossibility by constructing a quadruple $(\overline{\mathcal{N}}, r^*, \theta^*, \psi^*)$ with $R^\infty(\psi^*) < R^\infty(\psi)$.

Let $M_n$ be any countable model that contains $C$ and is c-isolated over $C$. Increase $B_\zeta$ for $\zeta \in I^+$ and extend the maps $f_\zeta$ so that $f_\zeta : M_n \to B_\zeta$. In particular, $B_\zeta$ is now a model for all $\zeta \in I$. Let $N_{\text{no }\eta}$ be prime over $B^*_{\text{no }\eta}$. Under the conditions we are working we can find $\hat{c} \in N_{\text{no }\eta}$ which realizes $r$ (it is a witness for $\chi$ with suitably chosen parameters). Let $\hat{C} = \text{acl}(M_{n-1}\hat{c})$.

We pause to introduce some more notation that will be used below. Fix an enumeration $\zeta_1, \ldots, \zeta_k$ of $I_\eta^+$.

- For any $q \in S(C)$ extending $p(\eta)$, let $q^*(x, z) = \text{tp}(ac/M_{n-1})$ for any realization $a$ of $q$.

- Let $q_i^* = f_{\zeta_i}(q^*)$. Formally, $q_i^*$ is over $f_{\zeta_i}(M_{n-1})$, but in places we identify it with its nonforking extension to $N_{\text{no }\eta}$.

- For $1 \le i \le k$, let $c_i = f_{\zeta_i}(c)$ and let $C_i = \text{acl}(N_{\text{no }\eta}c_i)$.

- If $q_1, \ldots, q_k \in S(C)$ each extend $p(\eta)$, an $\eta$-*sequence in* $q_1, \ldots, q_k$ is a sequence $\langle a_1, \ldots, a_k \rangle$, where $a_i$ realizes $f_{\zeta_i}(q_i)$. Note that since each $q_i$ is c-isolated, it follows from Lemma B.4(1) that $\{a_1c_1, \ldots, a_kc_k\}$ is independent over $N_{\text{no }\eta}$. In particular, $\text{tp}(a_1c_1, \ldots, a_kc_k/N_{\text{no }\eta})$ depends only on $q_1 \ldots, q_k$.

- If $q \in S(C)$ extends $p(\eta)$, then an $\eta$-*sequence in* $q$ is an $\eta$-sequence $\langle a_1, \ldots, a_k \rangle$, where every $a_i$ realizes $f_{\zeta_i}(q)$.

By absorbing parameters from $N_{\text{no }\eta}$, $\varphi$ has the form $\varphi(x, y_1z_1, \ldots, y_kz_k)$. Our hypotheses imply that

$$(1) \qquad \varphi(x, a_1c_1, \ldots, a_kc_k) \text{ is consistent and implies } q^*(x, c^*)$$

for any $q \in S(C)$ extending $p(\eta)$ and any $\eta$-sequence $\langle a_1, \ldots, a_k \rangle$ in $q$.

Fix a type $q \in S(C)$ extending $p(\eta)$ and an $\eta$-sequence $\langle a_1, \ldots, a_k \rangle$ in $q$. Let $b$ be any realization of $\varphi(x, a_1c_1, \ldots, a_kc_k)$. Since $\theta(b, \hat{c})$, we can choose $l$ least such that there is an $\hat{C}$-definable function $f(\overline{z})$ and a sequence $\overline{d}$ of realizations of $\psi$ such that $b = f(\overline{d})$. By Lemma B.4(2) there is a model $\hat{N} \supseteq \hat{C}b$ such that $\hat{N}$ is dominated by $\hat{C}$ over $M_{n-1}$. Choose a sequence $\overline{d}'$ of realizations of $\psi$ from $\hat{N}$ such that $b = f(\overline{d}')$. Since $\theta(b, \hat{c})$ holds, $\hat{C}b$ is dominated by $\hat{c}$ over $M_{n-1}$, hence every $d_j' \in \bar{d}'$ satisfies $d_j' \underset{M_{n-1}}{\downarrow} \hat{c}$. Since $\hat{C} \subseteq N_{\text{no }\eta}$, it follows that there is $\psi^* \in L(N_{\text{no }\eta})$ such that $\psi^*(d_j')$ holds for all $j$, but $R^\infty(\psi^*) < R^\infty(\psi)$. In particular, $b \in \text{dcl}(N_{\text{no }\eta} \cup \psi^*(\mathfrak{C}))$. As well, by shrinking $\varphi$ as needed we may assume that every realization of $\varphi(x, a_1 z_1, \ldots, a_k z_k)$ is in $\operatorname{dcl}(N_{\text{no }\eta} \cup \psi^*(\mathfrak{C}))$ for every $\eta$-sequence $\langle a_1, \ldots, a_k \rangle$ in any sequence of types $q_1, \ldots, q_k$ extending $p(\eta)$.

For each $1 \leq i \leq k$, let $t_i$ denote the sequence of variables $y_i z_i$ and define an equivalence relation $E_i(y, y') \in L(C_i)$ as follows. $E_1(y, y')$ is given by

$$d_{q_2^*} t_2 \ldots d_{q_k^*} t_k \forall x(\varphi(x, yc_1, t_2, \ldots, t_k) \leftrightarrow \varphi(x, y'c_1, t_2, \ldots, t_k))$$

and the others are defined analogously. For example, $E_k(y, y')$ is given by

$$d_{q_1^*} t_1 \ldots d_{q_{k-1}^*} t_{k-1} \forall x(\varphi(x, t_1, \ldots, t_{k-1}, yc_k) \leftrightarrow \varphi(x, t_1, \ldots, t_{k-1}, y'c_k)).$$

Ostensibly, each of the formulas $E_i(y, y')$ depend on our choice of $q_1, \ldots, q_k$. However, by iterating Lemma B.4(3) there is an $L(C)$-formula $\gamma(x) \in q$ extending $p(\eta)$ such that for each $i$, $E_i(y, y')$ does not depend on our choice of $q_1, \ldots, q_k$ so long as each $q_j$ extends $\gamma$. So fix such a $\gamma$.

We claim that there is at least one $i$ and a consistent $L(C)$-formula $\gamma' \vdash \gamma$ such that the $L(C_i)$-formula

$$\alpha(y) := \exists x(f_{\zeta_i}(\gamma')(x) \wedge x/E_i = y)$$

has no isolated extensions over $C_i$. Indeed, if this were not the case, then there would be a consistent $L(C)$-formula $\gamma^* \vdash \gamma$ such that

$$\alpha_i(y) := \exists x(f_{\zeta_i}(\gamma^*)(x) \wedge x/E_i = y)$$

is isolated for all $1 \leq i \leq k$. In particular, there would be distinct types $p, q \in S(C)$ extending $\gamma^*$ and $\eta$-sequences $\langle d_1, \ldots, d_k \rangle$, $\langle e_1, \ldots, e_k \rangle$ of $p$ and $q$ respectively such that $E_i(d_i, e_i)$ holds for each $i$. But then it would follow that

$$\forall x(\varphi(x, d_1 c_1, \ldots, d_k c_k) \leftrightarrow \varphi(x, e_1 c_1, \ldots, e_k c_k)).$$

However, this would contradict (1) above.

Fix such an $i$ and a $\gamma' \vdash \gamma$. Without loss, we may take $i = 1$ and we may identify $\overline{\mathcal{M}}$ with its image under $f_{\zeta_1}$. Let $C^* = C_1$ and let $r^*$ be the nonforking extension of $r$ to $C^*$. After this identification, we have a model $N_{\text{no }\eta} \supseteq \overline{\mathcal{M}}$, a formula $\psi^*$ over $N_{\text{no }\eta}$ with $R^\infty(\psi^*) < R^\infty(\psi)$, and a formula

$$\theta^*(y) := \exists x(\gamma'(x) \wedge x/E_1 = y)$$

over $C^*$ with no isolated extensions over $C^*$.

We can construe $N_{\text{no }\eta}$ to be the universe of an $n$-chain $\overline{\mathcal{N}} = \langle N_j : j < n \rangle$ that freely extends $\overline{\mathcal{M}}$ by choosing the models $N_j$ to satisfy

1. $N_0 \subseteq N_1 \subseteq \ldots \subseteq N_{\text{no }\eta}$ and

2. $N_j$ is prime over $\bigcup \{B_\zeta : \zeta \in I_{\text{no }\eta}, \zeta \not\geq \zeta_1|_{j+1}\}$.

In order to demonstrate that the quadruple $(\overline{\mathcal{N}}, r^*, \theta^*, \psi^*)$ contradicts the minimality of $(\overline{\mathcal{M}}, r, \theta, \psi)$, it remains to show that every realization of $\theta^*$ is $\psi^*$-definable over $C^*$. Choose a realization $a^*$ of $\theta^*$ and choose a realization $a$ of $\gamma'$ with $a/E_1 = a^*$. Let $q = tp(a/C)$. Let $s = \text{tp}(a_2 c_2, \ldots, a_k c_k / N_{\text{no } \eta})$, where $\langle a_1, \ldots, a_k \rangle$ is any $\eta$-sequence in $q$. Let $\langle \overline{d}_j : j \in \omega \rangle$ be a Morley sequence in $s$ over $C^*$ and let $D = \bigcup \{\overline{d}_j : j \in \omega\}$. To establish $\psi^*$-definabilty over $C^*$ it suffices by the note following Definition 3.13 to show that any automorphism $\sigma$ of $\mathfrak{C}$ that fixes $C^* \cup D \cup \psi^*(\mathfrak{C})$ pointwise also fixes $a^*$. Fix such a $\sigma$ and choose $j$ such that $\overline{d}_j \underset{C^*}{\downarrow} a\sigma(a)$. However, since any realization of $\varphi(x, ac_1, \overline{d}_j)$ is in $\text{dcl}(N_{\text{no } \eta} \cup \psi^*(\mathfrak{C}))$, it is fixed by $\sigma$. Hence $E_1(a, \sigma(a))$ holds, so $\sigma(a^*) = a^*$. Thus, the quadruple $(\overline{\mathcal{N}}, r^*, \theta^*, \psi^*)$ contradicts the minimality of $(\overline{\mathcal{M}}, r, \theta, \psi)$, so the three cases in the proof of the density of $D$ are exhaustive. So $D$ is dense in $\mathcal{P}$ and our proof is complete. □

We now show that if a special type fails to have a prime model over a chain, then the set of leaves constructed in Proposition 3.14 is diverse. Further, if the type is trivial as well, then the set is diffuse.

COROLLARY 3.17. 1. *If there is a chain $\overline{\mathcal{M}}$ of length $n$ and a special type $r \in R(\overline{\mathcal{M}})$, with no prime model over $M_{n-1}c$ for a realization $c$ of $r$, then there is a free extension $\overline{\mathcal{M}}'$ of $\overline{\mathcal{M}}$ with a diverse family of $\text{Leaves}(\overline{\mathcal{M}}')$ of size continuum.*

2. *If there is a chain $\overline{\mathcal{M}}$ of length $n$ and a trivial, special type $r \in R(\overline{\mathcal{M}})$, with no prime model over $M_{n-1}c$ for a realization $c$ of $r$, then there is a free extension $\overline{\mathcal{M}}'$ of $\overline{\mathcal{M}}$ with a diffuse family of $\text{Leaves}(\overline{\mathcal{M}}')$ of size continuum.*

*Proof.* (1) Fix any special type $r \in R(\overline{\mathcal{M}})$ such that there is no prime model over $M_{n-1}c$ for some realization $c$ of $r$. Choose the free extension $\overline{\mathcal{M}}'$ of $\overline{\mathcal{M}}$ and the family $Y = \{N_\eta : \eta \in {}^\omega 2\}$ constructed in Proposition 3.14. We claim that $Y$ is diverse. Indeed, let $Z_1, Z_2$ be distinct subsets of ${}^\omega 2$. Without loss, there is $\eta \in Z_1 \smallsetminus Z_2$. As $s_\eta(x, z)$ is realized in $N_\eta$, it is surely realized in $N^*(Z_1) \oplus_{M_{n-2}} V$ for any $(\overline{\mathcal{M}}, Y)$-model $V$. However, it follows immediately from Proposition 3.14 that $s_\eta$ is omitted in $N^*(Z_2)$. So, if $n = 1$, there is nothing more to prove. On the other hand, if $n > 1$, then since the tree conjugates of $r$ are foreign to some $\varphi \in r$, Lemma 3.10 tells us that any realization $c^*$ of $r$ in $N^*(Z_2) \oplus_{M_{n-2}} V$ is contained in $N^*(Z_2) \oplus_{M_{n-2}} W$ where $W$ is the truncation of $V$. But Proposition 3.14 tells us that $s_\eta(x, c^*)$ is omitted in $N^*(Z_2) \oplus_{M_{n-2}} V$ for any such $c^*$. Hence $Y$ is diverse.

(2) Here, fix a trivial, special type $r \in R(\overline{\mathcal{M}})$ such that there is no prime model over $M_{n-1}c$ for some realization $c$ of $r$. Then as above, the family $Y = \{N_\eta : \eta \in {}^\omega 2\}$ constructed in Proposition 3.14 is diffuse. To see this, choose $\eta \neq \nu$ and some $(\overline{\mathcal{M}}, Y)$-model $V$. In order to show that $N_\eta \oplus_{M_{n-1}} V \not\cong_V N_\mu \oplus_{M_{n-1}} V$ it suffices to show that $s_\eta^*(x, z)$ is omitted in $N_\mu \oplus_{M_{n-1}} V$ where $s_\eta^*(x, z)$ is the nonforking extension of $s_\eta$ to $V$. By way of contradiction, assume that $d^*c^*$ realizes $s_\eta^*$ in $N_\mu \oplus_{M_{n-1}} V$. Then, by Lemma 3.11, $c^*$ would be in $N_\mu \oplus_{M_{n-1}} W$ where $W$ is the truncation of $V$, but as before Proposition 3.14 says that $s_\eta(x, c^*)$ is omitted in $N_\mu \oplus_{M_{n-1}} V$, which yields our contradiction. □

3.3. *The existence of strongly regular types.* We begin by describing a general forcing construction that allows us to find a large subset of any perfect set of types over any countable, algebraically closed set in a stable theory in which no type in the set is isolated over any independent set of realizations of the others.

LEMMA 3.18. *Assume that $A \subseteq B$, where $A$ is countable and algebraically closed, $B$ realizes countably many complete types over $A$, and $Q \subseteq S(A)$ is a perfect subset. Then there is a subset $R \subseteq Q$ of size $2^{\aleph_0}$ such that for any distinct $r_1, \ldots, r_n, s \in R$, if $\{c_1, \ldots, c_n\}$ is independent over $A$, does not fork with $B$ over $A$, and each $c_i$ realizes $r_i$, then $s$ is not isolated over $Bc_1 \ldots c_n$.*

*Proof.* Call a formula $\theta(x)$ $Q$-perfect if $\{q \in Q : \theta \in q\}$ is perfect and let $\mathcal{Q}$ be the set of $Q$-perfect formulas. Let $\mathcal{P}$ denote the set of all finite functions $p : {}^{<\omega}2 \to \mathcal{Q}$ such that $p(\mu) \vdash p(\eta)$ for all $\eta \vartriangleleft \mu$ in $\mathrm{dom}(p)$. For $p, q \in P$, we say $p \leq q$ if and only if $\mathrm{dom}(p) \supseteq \mathrm{dom}(q)$ and $p(\eta) \vdash q(\eta)$ for all $\eta \in \mathrm{dom}(q)$.

We describe two countable lists of dense subsets of $\mathcal{P}$ that we wish our filter $G$ to meet. The first list of dense sets will ensure that we build complete types (hence elements of $Q$) in the limit. For each $\theta \in \mathcal{Q}$ and $\eta \in {}^{<\omega}2$, let

$$D_{\theta,\eta} = \{p \in \mathcal{P} : p(\eta) \vdash \theta \text{ or } p(\eta) \vdash \neg\theta\}.$$

Since $p(\eta)$ is $Q$-perfect for every $p \in \mathcal{P}$, each $D_{\theta,\eta}$ is dense in $\mathcal{P}$.

*Definition* 3.19. Suppose that $p \in \mathcal{P}$ with $\mathrm{dom}(p) \subseteq {}^{\leq m}2$ and $\eta_1, \ldots, \eta_k$ are (not necessarily distinct) elements from ${}^m 2$. An $L(A)$-formula $\varphi(x_1, \ldots, x_k)$ is *decided positively by $p$ at $\eta_1, \ldots, \eta_k$* if $\varphi(c_1, \ldots, c_k)$ holds for all $A$-independent sets $\{c_1, \ldots, c_k\}$ satisfying $p(\eta_i) \in \mathrm{tp}(c_i/A) \in Q$ for each $i$. The formula $\varphi$ is *decided by $p$ at $\eta_1, \ldots, \eta_k$* if either $\varphi$ or $\neg\varphi$ is decided positively by $p$ at $\eta_1, \ldots, \eta_k$. □

*Fact.* For every $L(A)$-formula $\varphi(x_1, \ldots, x_k)$, every $p \in \mathcal{P}$ with $\mathrm{dom}(p) \subseteq {}^{\leq m}2$, and every collection $\eta_1, \ldots, \eta_k$ from ${}^m 2$, there is $q \leq p$ that decides $\varphi$ at $\eta_1, \ldots, \eta_k$ and satisfies $\mathrm{dom}(q) \subseteq {}^{\leq m}2$. In addition, $q$ may be chosen so that $q_\mu = p_\mu$ for all $\mu \in {}^m 2$ distinct from $\eta_1, \ldots, \eta_k$.

*Proof.* By induction on $k$. Assume that the $L(A)$-formula $\varphi(x_0, \ldots, x_k)$, a condition $p$ satisfying $\text{dom}(p) \subseteq {}^{\leq m}2$, and $\eta_0, \ldots, \eta_k$ from ${}^m 2$ are given. Choose a $Q$-perfect formula $\theta$ of least $R(-, \varphi, \aleph_0)$-rank subject to

$$\theta \vdash p(\eta_0) \quad \text{and} \quad Mult(\theta, \varphi(x_0; x_1, \ldots, x_k), \aleph_0) = 1.$$

(Since $A$ is algebraically closed and $p(\eta_0)$ is $Q$-perfect, such a $\theta$ exists.) It follows that

$$d_r x_0 \varphi \equiv d_s x_0 \varphi$$

for all $r, s \in Q$ with $\theta \in r \cap s$. Let $\psi(x_1, \ldots, x_k)$ be the $L(A)$-formula $d_r x_0 \varphi$ for any such $r$ and let $p' \leq p$ be such that $\text{dom}(p') \subseteq {}^{\leq m}2$ and $p'(\eta_0) \vdash \theta$. Then, by applying the inductive hypothesis to $\psi$, we get $q \leq p'$ such that $q$ decides $\psi$ at $\eta_1, \ldots, \eta_k$. It follows that $q$ decides $\varphi$ at $\eta_0, \ldots, \eta_k$. □

As notation, for every pair of $L(A)$-formulas $\delta(x, y_1, \ldots, y_k, z)$ and $\theta(x)$, and every type $r \in S(A)$ realized in $B$, let

$$\Gamma_{\delta, r, \theta}(y_1, \ldots, y_k) := d_r z \forall x (\delta(x, y_1, \ldots, y_k, z) \to \theta(x)).$$

For every $L(A)$-formula $\delta$ and every type $r \in S(A)$ realized in $B$, let

$D_{\delta,r} = \{p \in \mathcal{P} : \text{there is } m \text{ such that } \text{dom}(p) \subseteq {}^{\leq m}2$ and for all $\eta_1, \ldots, \eta_k, \mu$ from ${}^m 2$ satisfying $\mu \neq \eta_i$ for each $i$, there is an $L(A)$-formula $\theta$ such that

$$\neg [\Gamma_{\delta, r, \theta}(c_1, \ldots, c_k) \leftrightarrow \theta(d)]$$

for all $d$ realizing $p(\mu)$ and all $A$-independent $\{c_1, \ldots, c_k\}$ satisfying $p(\eta_i) \in \text{tp}(c_1/A) \in Q$ for each $i\}$.

*Claim.* Each $D_{\delta,r}$ is dense in $\mathcal{P}$.

*Proof.* Choose $p \in \mathcal{P}$ arbitrarily. Choose $m$ such that $\text{dom}(p) \subseteq {}^{\leq m}2$. It suffices to handle each choice of $\eta_1, \ldots, \eta_k, \mu \in {}^m 2$ separately. So fix some such $\eta_1, \ldots, \eta_k, \mu$. Since $Q$ is perfect, there is a formula $\theta$ such that both $p(\mu) \wedge \theta$ and $p(\mu) \wedge \neg \theta$ are $Q$-perfect. From the Fact above, there is $p' \leq p$ that decides $\Gamma_{\delta, r, \theta}$ at $\eta_1, \ldots, \eta_k$ and satisfies $\text{dom}(p') \subseteq {}^{\leq m}2$ and $p'(\mu) = p(\mu)$.

There are two cases. In either case, put $q(\eta_i) = p'(\eta_i)$ for each $i$ and put $q(\gamma) = p'(\gamma)$ for all $\gamma \notin \{\eta_1, \ldots, \eta_k, \mu\}$. Put

$$q(\mu) = \begin{cases} p(\mu) \wedge \neg \theta & \text{if } \Gamma_{\delta, r, \theta} \text{ is decided positively by } p \text{ at } \eta_1, \ldots, \eta_k \\ p(\mu) \wedge \theta & \text{otherwise} \end{cases}$$

Now $q \leq p$ and $q$ meets our requirement for $\eta_1, \ldots, \eta_k, \mu$. After repeating this process for all sequences $\eta_1, \ldots, \eta_k, \mu \in {}^m 2$ with $\mu \neq \eta_i$ for each $i$ we obtain some $q^* \leq p$ with $q^* \in D_{\delta,r}$. □

Now fix any filter $G$ such that $G \cap D \neq \emptyset$ for each of the dense sets $D$ mentioned above. For each $\eta \in {}^\omega 2$, let

$$p_\eta(x) = \{\theta : p(\nu) \vdash \theta \text{ for some } p \in G \text{ and some } \nu \lessdot \eta\}.$$

The first collection of dense sets is used to guarantee that each $p_\eta$ is a complete type. Also, as each $p(\eta)$ consists of $Q$-perfect formulas, each $p(\eta) \in Q$. As any type trivially isolates itself, the second collection of dense sets ensures us that $p(\eta) \neq p(\mu)$ for all $\eta \neq \mu$. Now suppose that $\eta_1, \ldots, \eta_k, \mu \in {}^\omega 2$ are distinct. Pick $\{c_i : 1 \leq i \leq k\}$ to be independent over $A$, nonforking with $B$ over $A$, where each $c_i$ realizes $p(\eta_i)$. Suppose that $\delta(x, c_1, \ldots, c_k, b)$ isolates a type in $Q$. We claim that this type is not $p(\mu)$. To see this, let $r = \text{tp}(b/A)$ and choose $m$ such that $\eta_1|_m, \ldots, \eta_k|_m, \mu|_m$ are distinct. Pick $p \in G \cap D_{\delta,r}$ and let $\theta$ be the formula witnessing that $p \in D_{\delta,r}$ for $\eta_1|_m, \ldots, \eta_k|_m, \mu|_m$. If $\theta$ is in the type isolated by $\delta(x, c_1, \ldots, c_k, b)$ then $\Gamma_{\delta,r,\theta}(c_1, \ldots, c_k)$ holds, hence $\Gamma_{\delta,r,\theta}$ is decided positively by $p$ at $\eta_1|_m, \ldots, \eta_k|_m$. So $p(\mu|_m) \vdash \neg\theta$. On the other hand, if $\theta$ is not in the type isolated by $\delta(x, c_1, \ldots, c_k, b)$ then $\neg\Gamma_{\delta,r,\theta}$ is decided positively by $p$ at $\eta_1|_m, \ldots, \eta_k|_m$, so $p(\mu|_m) \vdash \theta$. In either case the type isolated by $\delta(x, c_1, \ldots, c_n, b)$ is not $p_\mu$. □

LEMMA 3.20. *Suppose that we are given an na-chain $\overline{\mathcal{M}}$ of length $n$, an $L(M_{n-1})$-formula $\varphi$, and a perfect set $Q \subseteq R(\overline{\mathcal{M}})$ such that every $q \in Q$ is special via $\varphi$.*

1. *There is a diverse family $Y \subseteq \text{Leaves}(\overline{\mathcal{M}})$ of size continuum.*

2. *If, in addition, every $q \in Q$ is trivial, then there is a diffuse family $Y \subseteq \text{Leaves}(\overline{\mathcal{M}})$ of size continuum.*

*Proof.* (1) Suppose that $\overline{\mathcal{M}}$, $\varphi$, and $Q$ are given. Choose $R \subseteq Q$ as in Lemma 3.18, taking $M_{n-1}$ for $A$ and $W$ for $B$. (Note that the isomorphism type of $W$ does not depend on the choice of $Y$.) We may assume that for any $r \in R$, there is a prime model $N_r$ over $M_{n-1}c$ for any realization $c$ of $r$, else there would already be a diverse family of size continuum by Corollary 3.17(1). We claim that $Y = \{N_r : r \in R\}$ is diverse. To see this, suppose $Z_1, Z_2 \subseteq R$ with $s \in Z_1 \smallsetminus Z_2$. It suffices to show that $s$ is omitted in $N^*(Z_2) \oplus_{M_{n-2}} V$ for any $(\overline{\mathcal{M}}, Y)$-model $V$. When $n = 1$ this follows immediately from Lemma 3.18, while if $n > 1$, then Lemma 3.10 tells us that any potential realization of $s$ must lie in $N^*(Z_2) \oplus_{M_{n-2}} W_V$, which it does not.

(2) Now suppose that the types are trivial as well. Arguing as above, apply Lemma 3.18 to obtain a subset $R \subseteq Q$ for $M_{n-1}$ and $W$. Further, by Corollary 3.17(2) we may assume that there is a prime model $N_r$ over any realization of $r$. In this case, the family $Y = \{N_r : r \in R\}$ will be diffuse. To see this, choose $s \neq r$ and let $s^*$ be the nonforking extension of $s$ to $V$, some $(\overline{\mathcal{M}}, Y)$-model. We claim that $s^*$ is omitted in $N_r \oplus_{M_{n-1}} V$. If it were realized by an element $a$, then Lemma 3.11 would imply that $a \in N_r \oplus_{M_{n-1}} W_V$, which would contradict Lemma 3.18. □

PROPOSITION 3.21. *Fix a chain $\overline{\mathcal{M}}$ of length $n$.*

1. *If some free extension $\overline{\mathcal{M}}'$ of $\overline{\mathcal{M}}$ has no diverse family of $\mathrm{Leaves}(\overline{\mathcal{M}}')$ of size continuum, then for every $p \in R(\overline{\mathcal{M}})$ there is a strongly regular $q \in R(\overline{\mathcal{M}})$ nonorthogonal to $p$ with a prime model over $M_{n-1}c$ for any realization $c$ of $q$.*

2. *If some free extension $\overline{\mathcal{M}}'$ of $\overline{\mathcal{M}}$ has no diffuse family of $\mathrm{Leaves}(\overline{\mathcal{M}}')$ of size continuum, then for every trivial $p \in R(\overline{\mathcal{M}})$ there is a strongly regular $q \in R(\overline{\mathcal{M}})$ nonorthogonal to $p$ with a prime model over $M_{n-1}c$ for any realization $c$ of $q$.*

*Proof.* We first prove (2). Suppose that the $na$-chain $\overline{\mathcal{M}}'$ is a free extension of $\overline{\mathcal{M}}$ that has no diffuse family of $\mathrm{Leaves}(\overline{\mathcal{M}}')$ of size continuum. Fix a trivial type $p$ in $R(\overline{\mathcal{M}})$. Choose a formula $\varphi \in L(M_{n-1})$ of least $R^\infty$-rank among all special formulas nonorthogonal to $p$ and let

$$X = \{q \in S(M_{n-1}) : q \text{ trivial, weight-1}, \varphi \in q, \text{ and } q \perp M_{n-2}\}.$$

If $n = 1$ we delete this last condition. By Lemmas C.2 and C.5, $X$ is a nonempty $G_\delta$ subset of $S(M_{n-1})$.

*Claim.* There are only countably many nonorthogonality classes represented in $X$.

*Proof.* If not, then since nonorthogonality is a Borel equivalence relation, Lemma C.1(1) implies that $X$ would contain a pairwise orthogonal family $\{q_i : i \in 2^{\aleph_0}\}$. Let $Z$ be the set of nonforking extensions of each $q_i$ to $M'_{n-1}$. For each $i$, choose a regular type $r_i$ domination-equivalent to $q_i$. Since $M'_{n-1}$ is an $na$-substructure of the universe, there is a regular type $s_i$ over $M'_{n-1}$ nonorthogonal to $q_i$. Clearly, each $s_i \in R(\overline{\mathcal{M}}')$, so we have a pairwise orthogonal family of regular types in $R(\overline{\mathcal{M}}')$ of size continuum, hence there is a diffuse family of $\mathrm{Leaves}(\overline{\mathcal{M}}')$ by Lemma 3.6, which is a contradiction. □

Let $\{r_i : i \in j \leq \omega\}$ be a maximal pairwise orthogonal subset of $X$ in which every $r_i \perp p$ and let

$$X_p = \{q \in X : q \not\perp p\}.$$

Since $q \in X_p$ if and only if $q \in X$ and $q \perp r_i$ for all $i \leq j$, it follows from Lemma C.4 that $X_p$ is a $G_\delta$ subset of $S(M_{n-1})$ as well. Since $\varphi$ was chosen to be special of least $R^\infty$-rank, it follows from the proof of Proposition 8.3.5 of [15] that every element of $X_p$ is regular as well. So, if $X_p$ were uncountable, then again by Lemma C.1(1) it would contain a perfect subset. By taking nonforking extension of each of these to $M'_{n-1}$, we would obtain a perfect subset of $R(\overline{\mathcal{M}}')$, where each element is special via $\varphi$. But this contradicts Lemma 3.20(2). Thus, we may assume that $X_p$ is a nonempty, countable $G_\delta$ subset of $S(M_{n-1})$. Hence, by Lemma C.1(2) there is an $L(M_{n-1})$-formula $\psi \vdash \varphi$ isolating some $q \in X_p$. By Proposition D.15 of [13], this $q$ is a strongly regular type via $\psi$. As well, since $q$ is special, if there were no prime model over $M_{n-1}c$ for a realization $c$ of $q$, then Corollary 3.17(2) would give us a diffuse family of Leaves$(\overline{\mathcal{M}})$ of size continuum.

We now prove (1). Again, assume that $\overline{\mathcal{M}}'$ is a na-chain freely extending $\overline{\mathcal{M}}$ with no diverse subset of Leaves$(\overline{\mathcal{M}}')$ of size continuum. Fix a type $p$ in $R(\overline{\mathcal{M}})$. If $p$ is trivial, then as we will see in Lemma 4.2 that every diffuse family is diverse, we are done by (2). Thus, we may assume that $p$ is nontrivial. Hence, by Lemma 8.2.20 of [15], there is a type $p' \in R(\overline{\mathcal{M}})$ nonorthogonal to $p$ and a formula $\theta \in p'$ witnessing that $p'$ is special and, in addition, the $p$-weight of $\theta$ is 1, and $p$-weight is definable inside $\theta$. Without loss, assume that $p = p'$. Let

$$X = \{q \in S(M_{n-1}) : \theta \in q, w_p(q) = 1\}.$$

Since $w_p(\theta) = 1$ and $p$-weight is definable inside $\theta$, $X$ is a closed subset of $S(M_{n-1})$. Now fix a $q \in X$. Since $w_p(q) = 1$, $q \not\perp p$. As well, since $q$ is $p$-simple, $q$ is regular. Hence, if $n > 1$ then $q \perp M_{n-2}$. That is,

$$X = \{q \in R(\overline{\mathcal{M}}) : q \text{ regular}, \theta \in q, \text{ and } q \not\perp p\}$$

is closed, hence a $G_\delta$. The proof is now analogous to (1). If $X$ were uncountable, then by Lemma C.1(1) it would contain a perfect subset, so by taking nonforking extensions to $M'_{n-1}$ there would be a diverse family of Leaves$(\overline{\mathcal{M}}')$ of size continuum by Lemma 3.20. So, we may assume that $X$ is a nonempty, countable, $G_\delta$ subset of $S(M_{n-1})$ and hence there is an $L(M_{n-1})$-formula $\psi \vdash \theta$ isolating some $q \in X$. Again by Proposition D.15 of [13], this $q$ is a strongly regular type via $\psi$. As well, since $q$ is special, if there were no prime model over $M_{n-1}c$ for a realization $c$ of $q$, then Corollary 3.17(1) would yield a diverse family of Leaves$(\overline{\mathcal{M}}')$ of size continuum. $\square$

### 3.4. Structure theorems for locally t.t. theories.

In this subsection we analyze some of the positive consequences of a theory being locally t.t. over an $n$-chain $\overline{\mathcal{M}}$.

The following three lemmas are true for any superstable theory. There is nothing novel about their statements or proofs, as they are central in the analysis of models of a totally transcendental theory. They are given here simply to indicate that the assumption of being totally transcendental can be weakened to include our context.

LEMMA 3.22.  *Suppose that $p, q \in S(M)$ where $p$ is regular and $q$ is strongly regular. Then $p$ and $q$ are orthogonal if and only if they are almost orthogonal over $M$.*

*Proof.* Suppose that $q$ is strongly regular via $\psi$ and that $p \not\perp q$. Choose a finite $e$ and realizations $a$ of $p|Me$ and $b$ of $q|Me$ respectively such that $a \underset{M}{\not\downarrow} eb$. Choose a formula $\varphi(x; y, z)$ over $M$ such that $\varphi(a; e, b)$ holds and $\varphi(a; y, z)$ forks over $M$. Thus
$$\exists z[\varphi(a; e, z) \wedge \psi(z)].$$
Since $M$ is a model and $a \underset{M}{\downarrow} e$, there is $e' \in M$ such that $\exists z[\varphi(a; e', z) \wedge \psi(z)]$. Let $b'$ be any witness to this formula. As $b'$ forks with $a$ over $M$, $\text{tp}(b'/M) \not\perp p$, hence $\text{tp}(b'/M) = q$ since $q$ is strongly regular via $\psi$. Thus $p$ and $q$ are not almost orthogonal over $M$. $\square$

LEMMA 3.23.  *Suppose that $p, q \in S(M)$ are nonorthogonal, where $p$ is regular and $q$ is strongly regular. Then any model containing a realization of $p$ contains a realization of $q$.*

*Proof.* Let $a$ realize $p$ and let $N \supseteq Ma$. Again suppose that $q$ is strongly regular via $\psi$. Since $p$ and $q$ are not almost orthogonal over $M$, there is a formula $\varphi(x, y)$ over $M$ such that $\varphi(a, y)$ forks over $M$ and $\varphi(a, b)$ holds for some $b$ realizing $q$. So $\exists y[\varphi(a, y) \wedge \psi(y)]$ holds in the monster model, hence by elementarity, $N$ contains a witness $b'$ to $\varphi(a, y) \wedge \psi(y)$. As before, it follows from the strong regularity of $q$ and the fact that $\varphi(a, y)$ forks over $M$ that $\text{tp}(b'/M) = q$. $\square$

*Notation.* Suppose that $p \in S(M)$. We adopt the notation $M(p)$ for the prime model over $M$ and any realization of $p$, if this prime model exists.

LEMMA 3.24.  *Suppose $q \in S(M)$ is strongly regular and $M(q)$ exists. Then $M(q)$ is atomic over $M$ and any $a \in M(q)$ such that $\text{tp}(a/M)$ is regular.*

*Proof.* Suppose that $M(q)$ is atomic over $Mb$, where $b$ realizes $q$, and that $q$ is strongly regular via $\psi(y)$. Choose $a \in M(q) \setminus M$ arbitrarily such that $tp(a/M)$ is regular. It suffices to show that $\text{tp}(b/Ma)$ is isolated. Let $\theta(x, b)$ isolate $\text{tp}(a/Mb)$. Since $a \notin M$, it follows from the Open Mapping Theorem that $a$ forks with $b$ over $M$, so let $\varphi(a,y) \in \operatorname{tp}(b/Ma)$ fork over $M$. Then the formula
$$\alpha(a,y) := \theta(a,y) \wedge \varphi(a,y) \wedge \psi(y)$$
isolates $\operatorname{tp}(b/Ma)$. To see this, it is evident that $\alpha(a,b)$ holds. So suppose that $\alpha(a,c)$ holds for some element $c$. We will show that $\operatorname{tp}(b/Ma) = \operatorname{tp}(c/Ma)$. Since $\varphi(a,y)$ forks over $M$, $\operatorname{tp}(a/M)$ is regular and $q$ is strongly regular via $\psi$, it must be that $\operatorname{tp}(c/M) = q$. Hence, $\theta(x,c)$ isolates a complete type over $Mc$. That is, $\operatorname{tp}(ab/M) = \operatorname{tp}(ac/M)$, so $\operatorname{tp}(b/Ma) = \operatorname{tp}(c/Ma)$ as desired. □

When we combine the three lemmas above with the existence of prime models over independent trees of models, we obtain a reasonable structure theory for the class of models of models extending a fixed chain. We illustrate this by continuing our analogy with the analysis of totally transcendental theories. The proofs of the following two corollaries follow immediately from the lemmas above.

COROLLARY 3.25. *If $M$ is countable and $p, q \in S(M)$ are both strongly regular and are nonorthogonal, then if $M(q)$ exists then $M(p)$ exists as well. Further, $M(p)$ and $M(q)$ are isomorphic over $M$.*

*Definition* 3.26. $I$ is a strongly regular sequence over $\overline{\mathcal{M}}$ if $I$ is independent over $M$ and every $a \in I$ realizes a strongly regular type in $R(\overline{\mathcal{M}})$.

COROLLARY 3.27. *If $T$ is locally* t.t. *over a chain $\overline{\mathcal{M}}$ of length $n$ and $I$ is a strongly regular sequence over $\overline{\mathcal{M}}$, then there is a prime model $N$ over $M_{n-1} \cup I$. Further, if $p$ and $q$ are nonorthogonal, strongly regular types from $R(\overline{\mathcal{M}})$, then $\dim(p, N) = \dim(q, N)$.*

The following lemma will be used in conjunction with the previous corollary to obtain good upper bounds.

LEMMA 3.28. *Suppose that $\overline{\mathcal{M}}$ is a d-chain, where d is the depth of $T$. If $T$ is locally* t.t. *over $\overline{\mathcal{M}}$, then any model $N \supseteq M_{d-1}$ with $N/M_{d-1} \perp M_{d-2}$ (when $d > 1$) is prime over $M_{d-1}$ and any maximal strongly regular sequence over $M_{d-1}$.*

*Proof.* Let $I$ be any maximal strongly regular sequence over $\overline{\mathcal{M}}$ inside $N$ and, from the previous corollary, let $N' \subseteq N$ be prime over $M_{d-1} \cup I$. Suppose, by way of contradiction, that $N' \neq N$. Choose $a \in N \setminus N'$ such that $p = \operatorname{tp}(a/N')$ is regular. Since $T$ has depth $d$, $p$ is not orthogonal to $M_{d-1}$, so choose a regular type $q$ over $M_{d-1}$ nonorthogonal to $p$. Since $T$ is locally t.t. over $\overline{\mathcal{M}}$, we may assume that $q$ is strongly regular. Let $q'$ be the nonforking extension of $q$ to $N'$. It follows from Lemma 3.23 that $q'$ is realized in $N$, but this contradicts the maximality of $I$. □

We conclude this subsection by giving a general lemma that will be used to aid our counting in some cases when the depth of the theory is small (typically $d = 1$ or $2$). As notation, note that any countable model $M$ can be viewed as a chain of length 1. In this case, we let $R(M)$ be the set of regular types over $M$ and we say that $T$ is locally t.t. over $M$ if it is locally t.t. over the associated chain of length 1.

LEMMA 3.29. *Suppose that $T$ is countable and superstable, and that $T$ is locally* t.t. *over any countable model $M$, and that $R(M)$ contains only countably many nonorthogonality classes. Then $T$ is $\omega$-stable.*

*Proof.* To show that $T$ is $\omega$-stable, for every countable model $M_0$, we will find a countable model $M_\omega$ realizing every type over $M_0$. So fix a countable model $M_0$. We construct $M_\omega$ to be a union of a chain of models $M_n$, where $M_0$ is given and for each $n$, $M_{n+1}$ is chosen to be prime over $M_n \cup I_n$, where $I_n$ is a strongly regular sequence over $M_n$ where each nonorthogonality class of $R(M_n)$ has dimension $\aleph_0$. (Such a sequence $I_n$ exists since $T$ is locally t.t. over $M_n$ and is countable since $R(M_n)$ is.)

We claim that $M_\omega$ realizes every type over $M_0$; in fact, it realizes every type over $M_n$ for all $n$. If not, then choose $p \in S(M_n)$ of least rank such that $p$ is omitted in $M_\omega$. Let $a$ be a realization of $p$ and let $N \supseteq M_n a$ be any countable model that is dominated by $a$ over $M_n$. (For instance, $N$ could be taken to be $l$-isolated over $M_n a$.) Let $q \in S(M_n)$ be any regular type realized in $N \setminus M_n$. Since $T$ is locally t.t. over $M_n$ there is a strongly regular type $r \in R(M_n)$ nonorthogonal to $q$, which by Lemma 3.23 is also realized in $N$, say by $b$. Since $N$ is dominated by $a$ over $M_n$, $a$ forks with $b$ over $M_n$, so $R^\infty(a/M_n b) < R^\infty(p)$. However, we can easily assume that $b \in M_{n+1}$ and so $R^\infty(a/M_{n+1}) < R^\infty(p)$ and certainly $tp(a/M_{n+1})$ is not realized in $M_\omega$ which is a contradiction. □

## 4. Unique decompositions and iteration

Before we can state one of the key theorems of this section, we introduce a very useful model $U$. Fix an $na$-chain $\overline{\mathcal{M}}$ of length $n$ and a set $Y$ of leaves of $\overline{\mathcal{M}}$. For every uncountable cardinal $\lambda$, we define an $(\overline{\mathcal{M}}, Y)$-model $U_\lambda$ which can be thought of as a '$\lambda$-saturated $(\overline{\mathcal{M}}, Y)$-model' of size $\lambda + 2^{\aleph_0}$. That is, we take $U_\lambda$ to be prime over an $(\overline{\mathcal{M}}, Y)$-tree $\mathcal{M} = \langle M_\zeta : \zeta \in I \rangle$. We require that $I$ is $\lambda$-branching above every $\eta \in I$ with $lg(\eta) < n$, i.e., $\{\nu : \nu^- = \eta\}$ has cardinality $\lambda$. Moreover, if $\langle f_\zeta : \zeta \in I \rangle$ is the family of elementary maps compatible with $\mathcal{M}$ demonstrating that it is an $(\overline{\mathcal{M}}, Y)$-tree, then for every $\zeta \in I$ with $lg(\zeta) = n$ and every $\bar{N} \in Y$, there are $\lambda$-many $\nu \in I$ with $\nu^- = \zeta$ such that $f_\nu$ maps onto $N_n$. For brevity, we denote $U_{\aleph_1}$ by $U$. The importance of this $(\overline{\mathcal{M}}, Y)$-model is given by the following theorem, whose proof follows from the main theorem of [7].

THEOREM 4.1 (Unique Decomposition Theorem). *If* $Q = \bigoplus_{M_k} \{P_i : i \in I\}$ *and* $Q' = \bigoplus_{M_k} \{P'_i : i \in I'\}$ *where each* $P_i, P'_i$ *is an* $(\overline{\mathcal{M}}, Y)$-*model with weight one over* $M_k$ *and* $Q \oplus_{M_k} V \cong_V Q' \oplus_{M_k} V$ *for some* $(\overline{\mathcal{M}}, Y)$-*model* $V$ *then there is a bijection* $f : I \to I'$ *such that* $P_i \oplus_{M_k} U \cong_U P'_{f(i)} \oplus_{M_k} U$ *for all* $i \in I$.

*Proof.* There is no loss in assuming that in fact $N = Q \oplus_{M_k} V = Q' \oplus_{M_k} V$ so in the terminology of [7] $\{P_i : i \in I\}$ and $\{P'_i : i \in I'\}$ are sets of independent $M_k$-components of $N$. There is a bijection $f : I \to I'$ so that for any $i \in I$, $\{P_i\} \cup \{P'_i : i \in I'\} \setminus \{P'_{f(i)}\}$ is a set of independent $M_k$-components of $N$. So if we fix $i \in I$ and let $\hat{V} = (\bigoplus_{M_k} \{P'_i : i \in I'\} \setminus \{P'_{f(i)}\}) \oplus_{M_k} V$ then $N = P_i \oplus_{M_k} \hat{V} = P'_{f(i)} \oplus_{M_k} \hat{V}$. Now choose $\lambda$ large enough so that $\hat{V}$ embeds into $U_\lambda$. By freely joining $U_\lambda$ with $N$ over $\hat{V}$ we obtain $P_i \oplus_{M_k} U_\lambda = P'_{f(i)} \oplus_{M_k} U_\lambda$. As well, note that $U \subseteq^{\mathcal{P}}_{\aleph_1} U_\lambda$ in the notation of [7]. Thus, by Theorem 2.9 of [7], $P_i \oplus_{M_k} U \cong_U P'_{f(i)} \oplus_{M_k} U$ which is the conclusion of the theorem. □

Our first application of this theorem gives the promised implication between the two central notions of this section.

LEMMA 4.2. *If* $Y$ *is diffuse then* $Y$ *is diverse.*

*Proof.* Assume that $Y$ is not diverse, i.e., there are distinct $Z_1, Z_2 \subseteq Y$ with $N^*(Z_1) \oplus_{M_{n-2}} V \cong N^*(Z_2) \oplus_{M_{n-2}} V$ over $M_{n-1}$ and $V$ for some $(\overline{\mathcal{M}}, Y)$-model $V$. Let $V'$ be the $(\overline{\mathcal{M}}, Y)$-model formed by taking the prime model over $M_{n-1}$ and $V$. We can assume that $N^*(Z_1) \oplus_{M_{n-1}} V' \cong_{V'} N^*(Z_2) \oplus_{M_{n-1}} V'$. Since $Z_1 \neq Z_2$ it follows from the Unique Decomposition Theorem that there are distinct $N, N' \in Y$ such that $N \oplus_{M_{n-1}} U \cong_U N' \oplus_{M_{n-1}} U$, which implies that $Y$ is not diffuse. □

The following definition captures a notion of homogeneity that our models will possess.

*Definition* 4.3. If $W \subseteq Q$, a set $H$ *reflects* $Q$ *over* $W$ if $W \subseteq H \subseteq Q$ and for every countable $A \subseteq Q$ there is an automorphism $\sigma$ of $Q$ over $W$ with $\sigma(A) \subseteq H$.

The utility of this notion is given by the lemma below, which will be used in many different contexts in what follows.

LEMMA 4.4. *Suppose that $\mathcal{F} = \{P_i : i \in \theta\}$ is a family of models, pairwise nonisomorphic over $V \cup W$, where $V$ is countable, and there is a cardinal $\lambda \geq |W|$ such that $\lambda^{\aleph_0} < \theta$ and every $P_i$ has a reflecting subset over $W$ of size at most $\lambda$. Then there is a subfamily $\mathcal{F}' \subseteq \mathcal{F}$ of size $\theta$ such that $P \not\cong_W P'$ for distinct $P, P' \in \mathcal{F}'$.*

*Proof.* Say $i \sim j$ if and only if $P_i \cong_W P_j$. It suffices to show that every $\sim$-class has size at most $\lambda^{\aleph_0}$. By way of contradiction, assume that some $\sim$-class $C \subseteq \theta$ has size greater than $\lambda^{\aleph_0}$. Without loss, assume $0 \in C$. Let $H \subseteq P_0$ be reflecting over $W$. For each $i \in C$, let $f_i : P_i \to P_j$ be an isomorphism over $W$. By composing each $f_i$ with an automorphism over $W$ if needed, the fact that $H$ is reflecting allows us to assume that $f_i(V) \subseteq H$ for all $i \in C$. Since $|C| > \lambda^{\aleph_0}$, there are distinct $i, j \in C$ such that $f_i(v) = f_j(v)$ for each $v \in V$. Thus, $f_j^{-1} \circ f_i : P_i \to P_j$ is an isomorphism over $V \cup W$, which is a contradiction. □

The following definition is arranged to allow us us to 'step down' trees and to achieve a lower bound on $I(T, \aleph_\alpha)$. The reader should observe that the hypotheses on $\kappa, \lambda$ imply that the existence of an iterable family of models implies that each model has size greater than continuum.

*Definition* 4.5. A family $\mathcal{F} = \{P_i : i \in \theta\}$ of $(\overline{\mathcal{M}}, Y)$-models is *k-iterable* if there are cardinals $\kappa, \lambda \geq \aleph_0$ such that:

1. $\lambda^{\aleph_0} < \kappa$ and $\lambda^{\aleph_0} < \theta$;

2. each $P_i$ has size $\kappa$ and has a reflecting subset $H_i$ over $M_k$ of size $\lambda$; and

3. $P_i \oplus_{M_k} U \not\cong_U P_j \oplus_{M_k} U$ for all distinct $i, j \in \theta$.

LEMMA 4.6. *If there is a $k$-iterable family of $\theta$ models, each of size $\aleph_\alpha$, then $I(T, \aleph_\alpha) \geq \theta$.*

*Proof.* Let $\mathcal{F}$ be such a family. It follows immediately from the definition that the models are pairwise nonisomorphic over $M_k$ and there is a cardinal $\lambda$ satisfying $\lambda^{\aleph_0} < \aleph_\alpha$ such that every $P \in \mathcal{F}$ has a reflecting subset of size $\lambda$ over $M_k$. Thus, $I(T, \aleph_\alpha) \geq \theta$ follows immediately from Lemma 4.4 (taking $W = \emptyset$). □

The intuition is that if a family is $k$-iterable, then we can use the lemma below to 'step down' the tree $k$ times, roughly exponentiating the number of models at each step.

LEMMA 4.7. *If $\mathcal{F} = \{P_i : i \in \theta\}$ is a k-iterable family of models each of size $\aleph_\alpha$ such that $k > 0$, $\theta \leq \aleph_\alpha$, and $\theta^{\aleph_0} < |\alpha + \omega|^\theta$, then there is a $(k-1)$-iterable family of $|\alpha + \omega|^\theta$ models, each of size $\aleph_\alpha$.*

*Proof.* Fix a cardinal $\lambda \geq \aleph_0$ such that $\lambda^{\aleph_0} < \aleph_\alpha$ and, for each $i \in \theta$ choose a subset $H_i$ reflecting $P_i$ over $M_k$. By enlarging each $H_i$ slightly we may take it to be the prime model over a subtree of a decomposition tree of $P_i$.

We first claim that there is a subfamily $\mathcal{F}'$ of $\mathcal{F}$ of size $\theta$ such that

$$P \oplus_{M_{k-1}} U \not\cong_U P' \oplus_{M_{k-1}} U$$

for all distinct $P, P' \in \mathcal{F}'$. For, if this were not the case, there would be a set $X \subseteq \theta$ of size $\theta$ such that $P_i \oplus_{M_{k-1}} U \cong_U P_j \oplus_{M_{k-1}} U$ for all $i, j \in X$. It is easily checked that $H_i' = H_i \oplus_{M_{k-1}} U$ is of size at most $\lambda \cdot 2^{\aleph_0}$ and reflects $P_i \oplus_{M_{k-1}} U$ over $U$. So, by applying Lemma 4.4 to the subfamily there would be $i \neq j$ such that $P_i \oplus_{M_{k-1}} U \cong_{UM_k} P_j \oplus_{M_{k-1}} U$. If we let $V$ be the prime model over $M_k$ and $U$, we can arrange that

$$P_i \oplus_{M_k} V \cong_V P_j \oplus_{M_k} V.$$

But this would imply $P_i \oplus_{M_k} U \cong_U P_j \oplus_{M_k} U$ by the Unique Decomposition Theorem, which would contradict the $k$-iterability of $\mathcal{F}$. Thus, by reindexing we may assume that our original family $\mathcal{F}$ satisfies $P_i \oplus_{M_{k-1}} U \not\cong_U P_j \oplus_{M_{k-1}} U$ for all distinct $i, j \in \theta$.

Let $S$ be an independent (over $M_{k-1}$) family of $\aleph_\alpha$ copies of $P_i$ for each $i \in \theta$. For each $P \in S$ choose a submodel $H_P$ prime over a subtree of a decomposition tree of $P$ of size at most $\lambda$. Let $J$ denote the set of all cardinal-valued functions $f : \theta \to \aleph_\alpha^+$ such that $f(0) = \aleph_\alpha$. For each $f \in J$, let $S_f \subseteq S$ consist of $f(i)$ copies of each $P_i$, and let $S_f^* \subseteq S_f$ consist of $\min\{f(i), \aleph_0\}$ copies of each $P_i$. Let $Q_f = \bigoplus_{M_{k-1}} S_f$. Each $Q_f$ has size $\aleph_\alpha$. Let $K_f = \bigoplus_{M_{k-1}} \{H_P : P \in S_f^*\}$. As any collection of automorphisms of distinct components of $Q_f$ extends to an automorphism of $Q_f$ over $M_{k-1}$, it follows that each $K_f$ reflects $Q_f$ over $M_{k-1}$. As well, $|K_f| \leq \theta \cdot \lambda$ and it follows from our cardinality assumptions that both $\lambda^{\aleph_0}$ and $\theta^{\aleph_0}$ are strictly less than $|\alpha + \omega|^\theta$.

Finally, suppose that $Q_f \oplus_{M_{k-1}} U \cong_U Q_g \oplus_{M_{k-1}} U$ for some distinct $f, g \in J$. Then, by the Unique Decomposition Theorem there are distinct $i, j \in \theta$ such that $P_i \oplus_{M_{k-1}} U \cong_U P_j \oplus_{M_{k-1}} U$, which is contrary to our additional hypothesis on $\mathcal{F}$ mentioned above. □

By iterating Lemma 4.7 we obtain the following:

LEMMA 4.8. *If for some $m \geq 1$ there is an m-iterable family of $\theta$ models, each of size $\aleph_\alpha$, where $\theta^{\aleph_0} = \theta$, then $I(T, \aleph_\alpha) \geq \min\{2^{\aleph_\alpha}, \beth_{m-1}(|\alpha + \omega|^\theta)\}$.*

*Proof.* Define a sequence of cardinals $\theta_0 < \theta_1 < \ldots < \theta_m$ by letting $\theta_0 = \theta$ and $\theta_{k+1} = |\alpha + \omega|^{\theta_k}$. Note that $\theta_k^{\aleph_0} = \theta_k$ for each $k$. There are now two cases. First, if $\theta_{m-1} \leq \aleph_\alpha$, then by applying Lemma 4.7 $m$ times, one gets a family of $\theta_m = \beth_{m-1}(|\alpha + \omega|^\theta)$ 0-iterable models over $M_0$. Hence, $I(T, \aleph_\alpha) \geq \beth_{m-1}(|\alpha + \omega|^\theta)$.

For the second case, assume that $\aleph_\alpha < \theta_{m-1}$. Choose $k$ least such that $\aleph_\alpha < \theta_k$ ($k$ may be zero). By applying Lemma 4.7 $k$ times, we obtain an $(m-k)$-iterable family $\mathcal{G}$ of $\theta_k$ models, each of size $\aleph_\alpha$. Now if $2^{\aleph_\alpha} \leq \theta_k$, then $I(T, \aleph_\alpha) = 2^{\aleph_\alpha}$ by applying Lemma 4.6 to this family. However, if $\theta_k < 2^{\aleph_\alpha}$, then as $\aleph_\alpha^{\aleph_0} \leq \theta_k^{\aleph_0} = \theta_k < 2^{\aleph_\alpha}$, we can apply Lemma 4.7 one more time to a subfamily of $\mathcal{G}$ of size $\aleph_\alpha$. This produces an $(m-k-1)$-iterable family of size $2^{\aleph_\alpha}$, hence $I(T, \aleph_\alpha) = 2^{\aleph_\alpha}$ by Lemma 4.6. $\square$

## 5. The counting

In this section, we combine the dichotomies from Section 2 with the machinery in Section 3 to obtain lower and upper bounds in many situations. Our results will be strong enough to compute the uncountable spectra in all cases. As noted in the introduction, we need only concern ourselves with (countable) classifiable theories of finite depth. For notation, *assume that such a theory $T$ has depth $d$.* We recall the nomenclature of [8].

*Definition* 5.1.  Suppose $1 \leq n \leq d$.

- $TT(n)$ holds if $T$ is locally t.t. over every chain $\overline{\mathcal{M}}$ of length $n$.

- $TF(n)$ holds if $T$ admits a trivial failure over some chain of length $n$ (i.e., some trivial $p \in R(\overline{\mathcal{M}})$ is not t.t. over $\overline{\mathcal{M}}$ for some chain $\overline{\mathcal{M}}$ of length $n$.)

- $NTF(n)$ holds (read: nontrivial failure) if there is a chain $\overline{\mathcal{M}}$ of length $n$ and a nontrivial type $p \in R(\overline{\mathcal{M}})$ that is not totally transcendental over $\overline{\mathcal{M}}$.

- We write $\#RD(n) = 2^{\aleph_0}$ if there is some chain $\overline{\mathcal{M}}$ of length $n$ where $R(\overline{\mathcal{M}})$ contains a family of continuum pairwise nonorthogonal types.

- We write $\#RD(n) = \aleph_0$ if there is some chain $\overline{\mathcal{M}}$ of length $n$ for which $R(\overline{\mathcal{M}})$ contains infinitely many nonorthogonal types, yet there is no chain $\overline{\mathcal{N}}$ of length $n$ with continuum many nonorthogonality classes represented in $R(\overline{\mathcal{N}})$.

- We say $\#RD(n)$ is finite if only finitely many nonorthogonality classes are represented in $R(\overline{\mathcal{M}})$ for every chain $\overline{\mathcal{M}}$ of length $n$.

- We write $\#RD(n) = 1$ if all types in $R(\overline{\mathcal{M}})$ are nonorthogonal for all chains $\overline{\mathcal{M}}$ of length $n$.

The notation $\#RD$ stands for the number of relevant dimensions. Since for any $n$-chain $\overline{\mathcal{M}}$, $S(M_{n-1})$ is a Polish space and the relation of nonorthogonality is a Borel equivalence relation, the number of nonorthogonality classes represented in $R(\overline{\mathcal{M}})$ is either countable or of size $2^{\aleph_0}$; hence at least one of the last four conditions hold for each $n$.

In the first subsection we use the machinery established earlier to obtain a number of lower bounds. Then, in the subsections that follow, we put on more and more conditions on our theory to obtain better and better upper bounds.

5.1. *Lower bounds.* The goal of this subsection is to obtain good lower bounds from instances of 'nonstructure' of the theory. The first of these is due to Shelah and has been known for some time. It will be used in two places, both where the depth of the theory is quite low. For a proof, see either Theorem 1.20 of Chapter IX of [18] or Theorem C of [2].

LEMMA 5.2. *If $T$ is superstable but not $\omega$-stable then $I(T, \aleph_\alpha) \geq \min\{2^{\aleph_\alpha}, \beth_2\}$ for all $\alpha > 0$.*

The next result is the 'General lower bound' mentioned in the Introduction. It is proved by the method of quasi-isomorphisms. (See either Theorem 5.10(a) of [16] or the discussion on page 396 of [1] for a proof.)

LEMMA 5.3. *If $T$ is classifiable of depth $d > 1$ then*
$$I(T, \aleph_\alpha) \geq \min\{2^{\aleph_\alpha}, \beth_{d-2}(|\alpha + \omega|^{|\alpha+1|})\}$$
*for all $\alpha > 0$.*

For each of the next five lemmas and propositions, *assume that $\overline{\mathcal{M}}$ is an na-chain of length $n$.* We begin by considering the effect of a diffuse family of leaves.

LEMMA 5.4. *If there is a diffuse family of $\mathrm{Leaves}(\overline{\mathcal{M}})$ of size $\mu \geq \aleph_0$, then for any ordinal $\alpha$ satisfying $\aleph_\alpha > 2^{\aleph_0}$ and $\mu^{\aleph_0} < |\alpha + \omega|^\mu$, there is an $(n-1)$-iterable family of $|\alpha + \omega|^\mu$ models, each of size $\aleph_\alpha$.*

*Proof.* Let $\mathcal{F} = \{P_i : i \in 2^{\aleph_0}\}$ be diffuse. Let $J$ be the set of all cardinal-valued functions $f : 2^{\aleph_0} \to \aleph_\alpha^+$ such that $f(0) = \aleph_\alpha$. It follows immediately from the Unique Decomposition Theorem that the family $\mathcal{G} = \{Q_f : f \in J\}$, where
$$Q_f = \bigoplus_{M_{n-1}} \{P_i^{(f(i))} : i \in 2^{\aleph_0}\}$$
satisfies $Q_f \oplus_{M_{n-1}} U \not\cong_U Q_g \oplus_{M_{n-1}} U$ for distinct $f, g \in J$. As well, the substructure $H_f = \bigoplus_{M_{n-1}} \{P_i^{(\hat{f}(i))} : i \in 2^{\aleph_0}\}$, where $\hat{f}(i) = \min\{f(i), \aleph_0\}$, reflects $Q_f$ over $M_{n-1}$ and has size continuum. $\square$

The following lemma is routine.

LEMMA 5.5. *If there is a diverse family of* $\text{Leaves}(\overline{\mathcal{M}})$ *of size* $\mu \geq \aleph_0$, *then for any cardinal* $\kappa \geq \mu$ *there is a family of* $2^\mu$ *pairwise nonisomorphic models over* $M_{n-1}$, *each of size* $\kappa$.

PROPOSITION 5.6. *If there is a diffuse family of* $\text{Leaves}(\overline{\mathcal{M}})$ *of size continuum then* $I(T, \aleph_\alpha) \geq \min\{2^{\aleph_\alpha}, \beth_{n-1}(|\alpha+\omega|^{2^{\aleph_0}})\}$ *for all ordinals* $\alpha > 0$.

*Proof.* Fix $\alpha > 0$. There are three cases. First, assume $2^{\aleph_\alpha} = 2^{\aleph_0}$. The existence of a diffuse family of size continuum clearly implies that the theory $T$ is not $\aleph_0$-stable, hence $I(T, \aleph_\alpha) \geq 2^{\aleph_0} = 2^{\aleph_\alpha}$ by Lemma 5.2. Second, assume $\aleph_\alpha \leq 2^{\aleph_0} < 2^{\aleph_\alpha}$. Then by Lemma 5.5 there is a family of $2^{\aleph_\alpha}$ models over $M_{n-1}$, each of size $\aleph_\alpha$, that are pairwise nonisomorphic over $M_{n-1}$. Our cardinal assumptions imply that $\aleph_\alpha^{\aleph_0} < 2^{\aleph_\alpha}$, so $I(T, \aleph_\alpha) = 2^{\aleph_\alpha}$.

Finally, if $\aleph_\alpha > 2^{\aleph_0}$, then by Lemma 5.4 there is an $(n-1)$-iterable family of $|\alpha+\omega|^{2^{\aleph_0}}$ models, each of size $\aleph_\alpha$. If $n = 1$, then $I(T, \aleph_\alpha) \geq |\alpha+\omega|^{2^{\aleph_0}}$ by Lemma 4.6. However, if $n > 1$, then by Lemma 4.8

$$I(T, \aleph_\alpha) \geq \beth_{n-2}\left(|\alpha+\omega|^{|\alpha+\omega|^{2^{\aleph_0}}}\right) = \beth_{n-1}(|\alpha+\omega|^{2^{\aleph_0}}).$$

We now turn our attention to the existence of a diverse family of leaves.

LEMMA 5.7. *If there is a diverse family of size* $\mu \geq \aleph_0$ *and* $n > 1$, *then for any* $\aleph_\alpha > 2^{\aleph_0}$ *there is an* $(n-2)$-*iterable family of size* $\min\{2^{\aleph_\alpha}, |\alpha+\omega|^{2^\mu}\}$, *where each model has size* $\aleph_\alpha$.

*Proof.* Choose a cardinal $\theta \leq \min\{\aleph_\alpha, 2^\mu\}$ such that

$$|\alpha+\omega|^\theta = \min\{2^{\aleph_\alpha}, |\alpha+\omega|^{2^\mu}\}, \quad \theta^{\aleph_0} < \aleph_\alpha, \quad \text{and} \quad \theta^{\aleph_0} < |\alpha+\omega|^\theta.$$

(Such a $\theta$ can be chosen from $\{2^{\aleph_0}, \beth_2, \aleph_\alpha\}$.) It follows directly from the diversity of the family of leaves that there is a collection $\mathcal{F} = \{P_i : i \in \theta\}$ of models over $M_{n-1}$, each of size at most $\aleph_\alpha$, that satisfies $P_i \oplus_{M_{n-2}} U \not\cong_U P_j \oplus_{M_{n-2}} U$ for all distinct $i, j \in \theta$.

Let $S$ be an independent (over $M_{n-2}$) family of $\aleph_\alpha$ copies of $P_i$ for each $i \in \theta$. Let $J$ denote the set of all cardinal-valued functions $f : \theta \to \aleph_\alpha^+$ such that $f(0) = \aleph_\alpha$. For each $f \in J$, let $S_f \subseteq S$ consist of $f(i)$ copies of each $P_i$, and let $S_f^* \subseteq S_f$ consist of $\min\{f(i), \aleph_0\}$ copies of each $P_i$. Let $Q_f = \bigoplus_{M_{n-2}} S_f$. Each $Q_f$ has size $\aleph_\alpha$. Let $K_f = \bigoplus_{M_{n-2}} \{P : P \in S_f^*\}$. Again, as any collection of automorphisms of distinct components of $Q_f$ extends to an automorphism of $Q_f$ over $M_{n-2}$, it follows that each $K_f$ reflects $Q_f$ over $M_{n-2}$. Since $|K_f| = \theta$ our hypotheses on $\theta$ ensure that this family satisfies the cardinal hypotheses of iterability.

So, suppose that $Q_f \oplus_{M_{n-2}} U \cong_U Q_g \oplus_{M_{n-2}} U$ for some distinct $f, g \in J$. Then, by the Unique Decomposition Theorem there are distinct $i, j \in \theta$ such that $P_i \oplus_{M_{n-2}} U \cong_U P_j \oplus_{M_{n-2}} U$, which is impossible. □

PROPOSITION 5.8. *If there is a diverse family of* Leaves$(\overline{\mathcal{M}})$ *of size continuum then*

$$I(T, \aleph_\alpha) \geq \begin{cases} \min\{2^{\aleph_\alpha}, \beth_2\} & \text{if } n = 1 \\ \min\{2^{\aleph_\alpha}, \beth_{n-2}(|\alpha + \omega|^{\beth_2})\} & \text{if } n > 1 \end{cases}$$

*for all ordinals $\alpha > 0$.*

*Proof.* First, if $\aleph_\alpha \leq 2^{\aleph_0}$ then $I(T, \aleph_\alpha) = 2^{\aleph_\alpha}$ by splitting into the same two cases as in the proof of Proposition 5.6. As well, if $\aleph_\alpha \geq 2^{\aleph_0}$ and $n = 1$, then Lemma 5.5 implies $I(T, \aleph_\alpha) \geq \beth_2$.

So assume $\aleph_\alpha > 2^{\aleph_0}$ and $n > 1$. Then by Lemma 5.7 there is an $(n-2)$-iterable family of $\min\{2^{\aleph_\alpha}, |\alpha + \omega|^{\beth_2}\}$ models, each of size $\aleph_\alpha$. Thus, the proposition follows from Lemma 4.6 if $n = 2$, and from Lemma 4.8 if $n > 2$.

We finish this section by obtaining lower bounds that arise from having suitably large collections of pairwise orthogonal types.

PROPOSITION 5.9. *If $R(\overline{\mathcal{M}})$ contains an infinite family of pairwise orthogonal types for some na-chain $\overline{\mathcal{M}}$ of length $n$, then*

$$I(T, \aleph_\alpha) \geq \min\{2^{\aleph_\alpha}, \beth_{n-1}(|\alpha + \omega|^{\aleph_0})\}$$

*for all ordinals $\alpha > 0$.*

*Proof.* First, we claim that $I(T, \aleph_\alpha) \geq 2^{\aleph_0}$ for all $\alpha > 0$. To see this we split into cases: If $T$ is not $\omega$-stable, this follows from Lemma 5.2. If $T$ is $\omega$-stable it can be verified by examining the spectra given by Saffe [16].

Next, fix an *na*-chain $\overline{\mathcal{M}}$ of length $n$ for which $R(\overline{\mathcal{M}})$ contains infinitely many pairwise orthogonal types. By Lemma 3.6 there is a diffuse family of Leaves$(\overline{\mathcal{M}})$ of size $\aleph_0$. We split into two cases.

If $|\alpha + \omega|^{\aleph_0} > 2^{\aleph_0}$, then $\alpha > 2^{\aleph_0}$, hence $\aleph_\alpha > 2^{\aleph_0}$. Thus, Lemma 5.4 provides us with an $(n-1)$-iterable family of $|\alpha + \omega|^{\aleph_0}$ models over $M_{n-1}$, each of size $\aleph_\alpha$. So, we obtain our lower bound via Lemma 4.6 or Lemma 4.8.

On the other hand, assume that $|\alpha + \omega|^{\aleph_0} = 2^{\aleph_0}$. If $n = 1$, then $I(T, \aleph_\alpha) \geq 2^{\aleph_0}$ as noted above, so we assume that $n > 2$. There are now three subcases.

- If $\aleph_\alpha > 2^{\aleph_0}$, then there is an $(n-2)$-iterable family of size $\beth_3$, where each model has size $\aleph_\alpha$. As before, the bound follows from Lemmas 4.6 and Lemma 4.8.

- If $\aleph_\alpha \leq 2^{\aleph_0}$, but $2^{\aleph_\alpha} > 2^{\aleph_0}$, then the existence of a diverse family of size $\aleph_0$ implies that there is a family $\mathcal{F} = \{P_i : i \in \aleph_\alpha\}$ of countable models over $M_{n-1}$ such that $P_i \oplus_{M_{n-2}} U \not\cong_U P_j \oplus_{M_{n-2}} U$ for distinct $i, j \in \aleph_\alpha$. Varying the dimensions of each of these yields a family of $2^{\aleph_\alpha}$ models, pairwise nonisomorphic over $M_{n-2}$, each of of size $\aleph_\alpha$. However, as $\aleph_\alpha^{\aleph_0} = 2^{\aleph_0}$ in this case, we obtain $I(T, \aleph_\alpha) = 2^{\aleph_\alpha}$.

- Finally, if $2^{\aleph_\alpha} = 2^{\aleph_0}$ then $I(T, \aleph_\alpha) = 2^{\aleph_0}$ as in the first paragraph.  $\square$

Using the notation given at the beginning of this section, we can easily summarize our lower bound results.

THEOREM 5.10.  *Let $T$ be a countable, classifiable theory of depth $d \geq n$.*

1. *If $\#RD(n) = 2^{\aleph_0}$ then $I(T, \aleph_\alpha) \geq \min\{2^{\aleph_\alpha}, \beth_{n-1}(|\alpha + \omega|^{2^{\aleph_0}})\}$ for all ordinals $\alpha > 0$.*

2. *If $\#RD(n) = \aleph_0$ then $I(T, \aleph_\alpha) \geq \min\{2^{\aleph_\alpha}, \beth_{n-1}(|\alpha + \omega|^{\aleph_0})\}$ for all ordinals $\alpha > 0$.*

3. *If $TF(n)$ holds then $I(T, \aleph_\alpha) \geq \min\{2^{\aleph_\alpha}, \beth_{n-1}(|\alpha + \omega|^{2^{\aleph_0}})\}$ for all ordinals $\alpha > 0$.*

4. *If $NTF(n)$ holds then*

$$I(T, \aleph_\alpha) \geq \begin{cases} \min\{2^{\aleph_\alpha}, \beth_2\} & \text{if } n = 1 \\ \min\{2^{\aleph_\alpha}, \beth_{n-2}(|\alpha + \omega|^{\beth_2})\} & \text{if } n > 1 \end{cases}$$

*for all ordinals $\alpha > 0$.*

*Proof.* (1) Choose a chain $\overline{\mathcal{M}}$ of length $n$ such that $R(\overline{\mathcal{M}})$ contains a pairwise orthogonal family of types. By passing to a free extension of $\overline{\mathcal{M}}$, we may assume that $\overline{\mathcal{M}}$ is an $na$-chain. Thus, it follows from Lemma 3.6 that there is a diffuse family of $\text{Leaves}(\overline{\mathcal{M}})$ of size continuum, so the lower bound follows immediately from Proposition 5.6.

(2) If some chain $\overline{\mathcal{M}}$ has an infinite family of pairwise orthogonal types in $R(\overline{\mathcal{M}})$, then this property will be inherited by any $na$-chain that freely extends $\overline{\mathcal{M}}$, so the result follows immediately from Proposition 5.9.

(3) Suppose that some chain $\overline{\mathcal{M}}$ of length $n$ supports a trivial type $p \in R(\overline{\mathcal{M}})$ that is not totally transcendental above $\overline{\mathcal{M}}$. By Lemma 3.9 we may assume that $p$ is special. Since $p$ is not totally transcendental, it follows from Corollary 3.17(2) and Proposition 3.21(2) that there is a free extension $\overline{\mathcal{M}}'$ of $\overline{\mathcal{M}}$ with a diffuse family of $\text{Leaves}(\overline{\mathcal{M}}')$ of size continuum. So the bound follows from Proposition 5.6.

(4) This is analogous to (3). If $p \in R(\overline{\mathcal{M}})$ is not totally transcendental over $\overline{\mathcal{M}}$, then, again using Lemma 3.9, it follows from Corollary 3.17(1) and Proposition 3.21(1) that there is a free extension $\overline{\mathcal{M}}'$ of $\overline{\mathcal{M}}$ with a diverse family of Leaves($\overline{\mathcal{M}}'$) of size continuum. Thus, the lower bound is given by Proposition 5.8. □

5.2. *Upper bounds.* In this short subsection we state some definitions and recall a very useful theorem for obtaining upper bounds.

*Definition* 5.11. Suppose that $\overline{\mathcal{M}}$ is an $n$-chain.

- An $\overline{\mathcal{M}}$-*component* is any model $N$ such that $M_{n-1} \subseteq_{na} N$ and $wt(N/M_{n-1}) = 1$.

- We write $\#C_{\overline{\mathcal{M}}}^\alpha$ for the number of $\overline{\mathcal{M}}$-components of size at most $\aleph_\alpha$, up to isomorphism over $M_{n-1}$ and

- write $\#C_n^\alpha = \sup\{\#C_{\overline{\mathcal{M}}}^\alpha : \overline{\mathcal{M}} \text{ is a } (d-n)\text{-chain}\}$ *unless* $n = 0$ and $\#C_{\overline{\mathcal{M}}}^\alpha$ is finite for all $d$-chains. In this case, put $\#C_0^\alpha = 1$.

PROPOSITION 5.12. 1. $\#C_0^\alpha \leq 2^{\aleph_0}$ *for all* $\alpha$.

2. *If* $TT(d)$ *holds then* $\#C_0^\alpha = \#RD(d)$ *unless the latter is finite, in which case* $\#C_0^\alpha = 1$.

The second part of the above Proposition follows directly from the definition of $TT(d)$ and Corollary 3.25.

The following theorem is proved by inductively counting the number of components as we step down a decomposition tree.

THEOREM 5.13. *If $T$ is a countable, classifiable theory with finite depth $d$ then*
$$\#C_{i+1}^\alpha \leq |\alpha + \omega|^{\#C_i^\alpha} + 2^{\aleph_0}$$
*and*
$$I(T, \aleph_\alpha) \leq \beth_{d-i-1}(|\alpha + \omega|^{\#C_i^\alpha} + 2^{\aleph_0})$$

*Proof.* By downward induction on $i$, using the fact that every model is prime over a normal tree of countable, $na$-substructures. □

As a corollary, we obtain the naive upper bound given in the Introduction.

COROLLARY 5.14. *If $T$ is countable, classifiable and has finite depth $d$ then*

$$I(T, \aleph_\alpha) \leq \beth_{d-1}(|\alpha + \omega|^{2^{\aleph_0}})$$

Combining this upper bound with Theorem 5.10 yields the following spectrum:

COROLLARY 5.15. *If $TF(d)$ holds or $\#RD(d) = 2^{\aleph_0}$ then*

$$I(T, \aleph_\alpha) = \min\{2^{\aleph_\alpha}, \beth_{d-1}(|\alpha + \omega|^{2^{\aleph_0}})\}.$$

In what follows, the assumptions of the given subsection are indicated in the subsection heading. Under these assumptions, a general upper bound will be derived and at the end of each subsection, we will indicate under what conditions this upper bound is met.

5.3. $TF(d)$ *fails,* $\#RD(d) \leq \aleph_0$ *and* $d > 1$. The hypotheses imply that we have control of the trivial components at level $d$. In order to get a better upper bound, we recall a theorem of Shelah that allows us to control the number of nontrivial components as well. The following lemma is Lemma 4.5 of Chapter XIII of [18].

LEMMA 5.16. *Suppose that $T$ is countable and classifiable. If $M_0 \subseteq_{\aleph_1} M_i$ and $tp(M_i/M_0)$ has weight one and is nontrivial for $i = 1, 2$ and $tp(M_1/M_0)$ and $tp(M_2/M_0)$ are not orthogonal then $M_1 \cong M_2$ over $M_0$.*

Fix a $(d-1)$-chain $\overline{\mathcal{M}}$ and $N$, an $\overline{\mathcal{M}}$-component of size at most $\aleph_\alpha$. We wish to find a set of invariants which determines the isomorphism type of $N$ over $M_{d-2}$.

Choose $M_{d-1} \subseteq_{na} N$ so that $M_{d-1}$ is countable and properly contains $M_{d-2}$. Now let $I$ be a maximal $M_{d-1}$-independent collection of countable models contained in $N$, where each is a weight one, $na$-extension of $M_{d-1}$. Since $T$ has depth $d$, $N$ is prime over $I$. Since $TF(d)$ fails, each $N' \in I$ for which $tp(N'/M_{d-1})$ is not orthogonal to a trivial regular type, is actually determined, up to isomorphism over $M_{d-1}$ by this information alone (see Corollary 3.27). However, we do not have such control over the nontrivial types. In order to remedy this, we do the following.

Choose $I_0 \subseteq I$, $|I_0| \leq 2^{\aleph_0}$ such that $M'$, the prime model over $I_0$, is an $\aleph_1$-substructure of $N$. By Fact 5.16, the isomorphism types of the nontrivial components are determined, up to isomorphism over $M'$, by their nonorthogonality class. Since the nontrivial (and trivial) components present in $I$ are all based on $M_{d-1}$, there are $|\alpha+\omega|^{\#RD(d)}$ possibilities for the isomorphism type of $N$ over $M'$. We did have to fix a countable model $M_{d-1}$ over $M_{d-2}$ and a model $M'$ of size at most $2^{\aleph_0}$. Clearly, there are at most $\beth_2$ many possibilities for these choices, up to isomorphism over $M_{d-2}$.

In summary, if $TF(d)$ fails and $d > 1$ then $\#C_1^\alpha \leq |\alpha+\omega|^{\#RD(d)} + \beth_2$. Thus, by Theorem 5.13, we have the following:

$$I(T, \aleph_\alpha) \leq \beth_{d-2}(|\alpha+\omega|^{|\alpha+\omega|^{\#RD(d)}+\beth_2} + 2^{\aleph_0}) = \beth_{d-1}(|\alpha+\omega|^{\#RD(d)} + \beth_2).$$

Combining these upper bounds with the lower bounds of Theorem 5.10 yields the following corollaries.

COROLLARY 5.17. 1. If $TF(d)$ fails, $NTF(d)$ holds, $\#RD(d) = \aleph_0$, and $d > 1$ then
$$I(T, \aleph_\alpha) = \min\{2^{\aleph_\alpha}, \beth_{d-1}(|\alpha+\omega|^{\aleph_0} + \beth_2)\}.$$

2. If $TF(d)$ fails, $NTF(d)$ holds, $\#RD(d)$ is finite, and $d > 1$ then
$$I(T, \aleph_\alpha) = \min\{2^{\aleph_\alpha}, \beth_{d-1}(|\alpha+\omega| + \beth_2)\}.$$

*Proof.* For both cases, note that $I(T, \aleph_\alpha) \geq \min\{2^{\aleph_\alpha}, \beth_{d+1}\}$ since $NTF(d)$ holds. If $\#RD(d) = \aleph_0$ then the bound $I(T, \aleph_\alpha) \geq \min\{2^{\aleph_\alpha}, \beth_{d-1}(|\alpha+\omega|^{\aleph_0})\}$ follows from Theorem 5.10(2), Combining these lower bounds with the upper bound mentioned above yields the spectrum in (1).

On the other hand, if $\#RD(d)$ is finite, then by combining the lower bound of the preceding paragraph with the general lower bound (Lemma 5.3), we match the upper bound mentioned above. □

Thus, except for the case when $d = 1$, which we handle in Subsection 5.8, we have computed the spectra of all (classifiable of finite depth $d$) theories for which $TT(d)$ fails or $\#RD(d) = 2^{\aleph_0}$.

5.4. $TT(d)$ *holds and* $\#RD(d) \leq \aleph_0$. Since $TT(d)$ holds, the positive results from Subsection 3.4 come into play and yield substantially better upper bounds. In particular, Proposition 5.12 implies $\#C_0^\alpha = \aleph_0$ if $\#RD(d) = \aleph_0$ and 1 if $\#RD(d)$ is finite.

COROLLARY 5.18. If $TT(d)$ holds and $\#RD(d) = \aleph_0$ then
$$I(T, \aleph_\alpha) = \min\{2^{\aleph_\alpha}, \beth_{d-1}(|\alpha+\omega|^{\aleph_0})\}.$$

*Proof.* Since $\#C_0^\alpha = \aleph_0$, it follows from Theorem 5.13 that an upper bound is
$$\beth_{d-1}(|\alpha+\omega|^{\aleph_0} + 2^{\aleph_0}) = \beth_{d-1}(|\alpha+\omega|^{\aleph_0}).$$

However, the matching lower bound is immediate from Theorem 5.10(2). □

Suppose that $\#RD(d)$ is finite. Then Theorem 5.13 yields an upper bound of
$$\beth_{d-1}(|\alpha + \omega| + 2^{\aleph_0}).$$

Note that if $d > 1$ then we can obtain a matching lower bound if either $TT(d-1)$ fails or $\#RD(d-1) = 2^{\aleph_0}$.

COROLLARY 5.19. *If $TT(d)$ holds, $\#RD(d)$ is finite, $d > 1$, and either $TT(d-1)$ fails or $\#RD(d-1) = 2^{\aleph_0}$, then*
$$I(T, \aleph_\alpha) = \min\{2^{\aleph_\alpha}, \beth_{d-1}(|\alpha + \omega| + 2^{\aleph_0})\}.$$

*Proof.* If either $TF(d-1)$ holds or $\#RD(d-1) = 2^{\aleph_0}$, then the lower bound follows immediately from Theorem 5.10. If $NTF(d-1)$ holds, then the lower bound is obtained by combining the lower bound of Theorem 5.10(4) with the general lower bound of Lemma 5.3. □

5.5. *Obtaining na-inclusion.* Before continuing, we give a technical lemma and a construction that are needed to establish the more fussy upper bounds.

LEMMA 5.20. *Suppose that $\overline{\mathcal{M}}$ is a chain of length $n > 1$ and $TT(n)$ holds. If $M_{n-1} \subseteq N$ such that $tp(N/M_{n-1})$ is orthogonal to $M_{n-2}$, and every strongly regular type over $M_{n-1}$ which is realized in $N \setminus M_{n-1}$ has infinite dimension in $M_{n-1}$, then $M_{n-1} \subseteq_{na} N$.*

*Proof.* Suppose that $\varphi(x, \bar{a})$ has a solution in $N \setminus M_{n-1}$ where $\bar{a} \in M_{n-1}$. By relativizing the proof that some regular type is realized between any pair of models of a superstable theory, there is an element $d \in dcl(\varphi(N)) \setminus M_{n-1}$ such that $q = tp(d/M_{n-1})$ is regular. Since $TT(n)$ holds, there is a strongly regular type $p \in R(\overline{\mathcal{M}})$ nonorthogonal to $q$. Without loss, we may assume that $p$ is based on $\bar{a}$. So, by Lemma 3.23 there a realization $c$ of $p$ in $N \setminus M_{n-1}$ that depends on $d$ over $M_{n-1}$. Choose $b \in \varphi(N) \setminus M_{n-1}$ that depends on $c$ over $M_{n-1}$ and choose $\bar{e} \supseteq \bar{a}$ from $M_{n-1}$ such that $b$ and $c$ are dependent over $\bar{e}$. Let the formula $\chi(b, c, \bar{e})$ witness this dependence. Let $I$ be an infinite Morley sequence in the type of $p|\bar{a}$ inside $M_{n-1}$. Choose $c' \in I$ so that $c'$ and $\bar{e}$ are independent over $\bar{a}$. We have
$$\exists \chi(x, c, \bar{e}) \wedge \varphi(x, \bar{a})$$
which is also true when $c'$ replaces $c$. So pick $b'$ so that
$$\chi(b', c', \bar{e}) \wedge \varphi(b', \bar{a})$$
holds. Since $b'$ and $c'$ are dependent over $\bar{e}$, $b'$ cannot be in the algebraic closure of $\bar{a}$ and so we finish. □

We next describe a construction that will be used in the next subsections.

*Construction* 5.21. Fix a chain $\overline{\mathcal{M}}$ of length $(n-1)$ and an $\overline{\mathcal{M}}$-component $N$. Suppose that $n > 1$, $TT(n-1)$ and $TT(n)$ hold and $\#RD(n)$ is finite. Since $TT(n-1)$ holds, we can find $a \in N \setminus M_{n-2}$ such that $\operatorname{tp}(a/M_{n-2})$ is strongly regular. Let $N_0$ be contained in $N$ and prime over $M_{n-2}a$; again, $N_0$ exists because we are assuming $TT(n-1)$. Let $\overline{\mathcal{N}}_0$ be the chain of length $n$ obtained by concatenating $N_0$ to $\overline{\mathcal{M}}$. We will define, by induction on $i$, countable models $N_i \subseteq N$, numbers $n_i$ and strongly regular types $p_l^i \in R(N_i)$ for $l < n_i$ that are pairwise orthogonal.

So suppose that $N_i$ has been defined. Let $\overline{\mathcal{N}}_i$ be the chain of length $n$ obtained by concatenating $N_i$ to $\overline{\mathcal{M}}$. Let $n_i$ be the number of nonorthogonality classes in $R(\overline{\mathcal{N}}_i)$ that are realized in $N$ and that are orthogonal to $p_l^j$ for all $j < i$ and $l < n_j$. Since $TT(n)$ holds, we can find strongly regular representatives $p_l^i$ for $l < n_i$ of these classes.

Now, for $l < n_i$, let $I_l^i$ be a maximal Morley sequence in $N$ for $p_l^i$ if the dimension of this type in $N$ is countable and otherwise, let $I_l^i$ be any countable, infinite Morley sequence in $N$ for $p_l^i$. Let $I^i$ be the concatenation of the sequences $I_l^i$ and finally, let $N_{i+1}$ be prime over $N_i I^i$.

We next argue that this process must stop after finitely many steps. To see this, notice that if not, then if $N_\omega = \bigcup_i N_i$, $\overline{\mathcal{N}}_\omega$ together with all the $p_l^i$'s exemplifies that $\#RD(n)$ is not finite. What have we achieved? If $N_i$ is the last element of this chain then by Lemma 5.20, $N_i \subseteq_{na} N$. To see this, suppose that $b \in N \setminus N_i$ such that $q = tp(b/N_i)$ is strongly regular. Then for some $j$ and $l$, $q$ is not orthogonal to $p_l^j$ and the dimension of this latter type in $N_i$ is infinite.

5.6. $TT(d)$ and $TT(d-1)$ hold, $\#RD(d)$ is finite, $d > 1$, and $\#RD(d-1) \leq \aleph_0$. Our first goal is to obtain an upper bound for theories satisfying these hypotheses. This is accomplished by analyzing Construction 5.21 in detail.

Fix a chain $\overline{\mathcal{M}}$ of length $d-1$ and an $\overline{\mathcal{M}}$-component $N$ of size at most $\aleph_\alpha$. Pick $a \in N$ such that $tp(a/M_{d-2})$ is strongly regular.

Let $N'$ be the model described in Construction 5.21. $N' \subseteq_{na} N$ and so by Lemma 3.28, $N$ is prime over $N'$ and $I$, a strongly regular sequence over $N'$. By the construction of $N'$ there are $|\alpha + 1|^{\#RD(d)}$ possibilities for $I$ up to isomorphism over $N'$. Now the construction of $N'$ was accomplished in finitely many steps. We see that at each step $i$, there were at most countably many choices for $I^i$ and at the first stage, since $\#RD(d-1) \leq \aleph_0$, countably many choices for $N_0$ so there are at most countably many choices for $N'$.

In summary, there are at most $|\alpha + \omega|$ many isomorphism types of $N$ over $M_{d-2}$, i.e., $\#C_1^\alpha \leq |\alpha + \omega|$. Hence, under the assumptions of this subsection, we get an upper bound of

$$\beth_{d-2}(|\alpha + \omega|^{|\alpha+\omega|} + 2^{\aleph_0}) = \beth_{d-1}(|\alpha + \omega|).$$

Note that this upper bound matches the general lower bound (Lemma 5.3) whenever $\alpha$ is infinite. So, in the computation of the spectra that follow, we are only interested in computing $I(T, \aleph_\alpha)$ when $\alpha$ is finite. We can dispense with a number of cases at this point.

COROLLARY 5.22. *Suppose that $TT(d)$ and $TT(d-1)$ hold, $d > 1$ and $\#RD(d)$ is finite. If any of the following three conditions hold*

- $d > 1$ and $\#RD(d-1) = \aleph_0$;
- $d > 2$ and $TT(d-2)$ *fails; or*
- $d > 2$ and $\#RD(d-2) = 2^{\aleph_0}$

*then $I(T, \aleph_\alpha) = \min\{2^{\aleph_\alpha}, \beth_{d-1}(|\alpha + \omega|)\}$.*

*Proof.* All three cases follow immediately by combining lower bounds from Theorem 5.10 with the general lower bound (Lemma 5.3) and matching the upper bound mentioned above. □

To distinguish between the spectra $\beth_{d-1}(|\alpha+\omega|)$ and $\beth_{d-2}(|\alpha+\omega|^{|\alpha+1|})\}$, we need one further dichotomy.

*Definition* 5.23. Suppose that $T$ is a countable, classifiable theory with finite depth $d > 1$ and moreover, both $TT(d)$ and $TT(d-1)$ hold and $\#RD(d)$ is finite. We say that $T$ has the *final property* if for every chain $\overline{\mathcal{M}}$ of length $d$, there are only finitely many isomorphism types of models $N$ over $M_{d-2}$ of the form $Pr(M_{d-1} \cup J)$, where $J$ is a countable, strongly regular sequence from $RD(\overline{\mathcal{M}})$.

We remark that in the case where $T$ is $\omega$-stable (and satisfies the other properties) $T$ has the final property if and only if all types of depth $d - 2$ are abnormal (of Type $V$) in the sense of Baldwin [1].

COROLLARY 5.24. *Suppose that $d > 1$, $TT(d)$ and $TT(d-1)$ hold, and $\#RD(d)$ is finite. If $T$ does not have the final property, then*

$$I(T, \aleph_\alpha) = \min\{2^{\aleph_\alpha}, \beth_{d-1}(|\alpha + \omega|)\}.$$

*Proof.* First, if $T$ is totally transcendental, then this is proved in Saffe [16] or Baldwin [1], so we assume that $T$ is not totally transcendental. Hence, by Lemma 3.29, we may assume that $d > 2$. So, by the previous corollary, we may further assume that $\#RD(d-1)$ is finite, $TT(d-2)$ holds, and that $\#RD(d-2) \leq \aleph_0$. As well, if $\aleph_\alpha \leq 2^{\aleph_0}$, then we are done by Lemma 5.2, so choose $\alpha$ such that $\aleph_\alpha > 2^{\aleph_0}$. As noted above, we may assume $\alpha$ is finite (else the upper bound already matches the general lower bound). We need to show that

$$I(T, \aleph_\alpha) \geq \min\{2^{\aleph_\alpha}, \beth_{d-1}\}.$$

It is easily checked that if the final property fails for some chain, then it fails for some $na$-chain, so choose $\overline{\mathcal{M}}$ an $na$-chain of length $d$ for which the final property fails. Then, by using the method presented in the proof of Lemma 5.6 of Chapter XVIII of [1] one obtains a family $\{P_i : i \in 2^{\aleph_0}\}$ of countable models over $M_{d-2}$ such that $P_i \oplus_{M_{d-3}} U \not\cong_U P_j \oplus_{M_{d-3}} U$ for distinct $i, j$. Then, arguing as in the last two paragraphs of Lemma 5.7, we obtain a $(d-3)$-iterable family of models of size $\beth_2$, each of which is of size $\aleph_\alpha$. Hence, we obtain our lower bound from Lemmas 4.6 and 4.8. □

5.7. *The final case.* The assumptions of this case are too many to put in the subsection heading. We will assume that $d > 2$, $TT(d)$, $TT(d-1)$ and $TT(d-2)$ hold, $\#RD(d)$ and $\#RD(d-1)$ are finite, $\#RD(d-2) \leq \aleph_0$, and assume that the final property holds. The case when $d = 2$ is handled in the next subsection.

In order to compute the best general upper bound in this case, we will compute $\#C_2^\alpha$. Toward this end, fix a chain $\overline{\mathcal{M}}$ of length $d-2$ and an $\overline{\mathcal{M}}$-component $N$ of size at most $\aleph_\alpha$. Now fix $a \in N$ such that $tp(a/M_{d-3})$ is strongly regular and let $N_0$ be the prime model over $M_{d-3}a$. Since we are assuming that $TT(d-2)$ holds and $\#RD(d-2) \leq \aleph_0$, there are at most countably many choices for $N_0$ up to isomorphism over $M_{d-3}$. Now by using Construction 5.21, we can find $M_{d-2} \subseteq_{na} N$ and by using the argument from the previous subsection, there are at most countably many choices for $M_{d-2}$ over $N_0$. Let $\overline{\mathcal{M}}'$ be the chain of length $d-1$ formed by concatenating $M_{d-2}$ to $\overline{\mathcal{M}}$. Now in order to understand $N$ over $M_{d-2}$, it suffices to understand the isomorphism types of $\overline{\mathcal{M}}'$-components inside $N$. Let $N'$ be the model described in Construction 5.21, now working over $\overline{\mathcal{M}}'$. Again, $N' \subseteq_{na} N$ and $N$ is prime over $N'$ and a strongly regular sequence $I$ over $N'$. As before, there are only $|\alpha+1|^{\#RD(d)}$ possibilities for $I$ over $N'$. However, since $\#RD(d-1)$ is finite, there are only finitely many choices for the first model in the construction of $N'$. In addition, since the final property holds, then at each stage there are only finitely many choices for the $i^{th}$ model. And, as we noted, the construction of $N'$ is accomplished in finitely many steps and so by König's Lemma, there are only finitely many choices for $N'$ in all.

In summary, there are at most $|\alpha+\omega|$ choices for $M_{d-2}$ and $|\alpha+l|$ choices for $N$ over $M_{d-2}$. Hence, over $M_{d-2}$, there are at most $|\alpha+\omega|^{|\alpha+1|}$ many such $N$ up to isomorphism. We conclude then that

$$\#C_2^\alpha \leq |\alpha+\omega|^{|\alpha+1|} + \aleph_0 = |\alpha+\omega|^{|\alpha+1|}$$

and the upper bound in this case is

$$\beth_{d-3}(|\alpha+\omega|^{|\alpha+\omega|^{|\alpha+1|}} + 2^{\aleph_0}) = \beth_{d-2}(|\alpha+\omega|^{|\alpha+1|}).$$

As this agrees with the general lower bound, it is the spectrum in this case.

5.8. *The case of $d = 1$ or $2$.* In this section we would like to take care of the case when $d = 1$ and finish the previous subsection in the case when $d = 2$.

To improve the upper bounds in the case $d = 1$ we make several remarks. If $TF(1)$ fails and $\#RD(1) = 1$ (this case was handled in subsection 5.3), we actually obtain the upper bound $\beth_2$. This can be seen even in the proof presented in that subsection. After fixing a model of size $2^{\aleph_0}$, there is only one dimension possible. Of course, the spectrum $\min\{2^{\aleph_\alpha}, \beth_2\}$ is achieved when $T$ is unidimensional and not $\omega$-stable.

Next, we improve the upper bounds when $TT(1)$ holds and $\#RD(1)$ is finite. It follows directly from Lemma 3.29 that such theories are totally transcendental, so we could quote Saffe [16]. However, we include a brief discussion for completeness. The upper bound then, from subsection 5.3 is

$$|\alpha + \omega|.$$

If $\#RD(1)$ is not 1 then in [11] and [12], Lachlan shows that this spectrum is correct under these assumptions unless $T$ is $\omega$-categorical. In those papers, he shows that if $T$ is $\omega$-categorical then for some number $m$ and some $G \leq \text{Sym}(m)$,

$$I(T, \aleph_\alpha) = |(\alpha+1)^m/G| - |(\alpha)^m/G|$$

Of course, the case when $TT(1)$ holds and $\#RD(1) = 1$ is when $T$ is $\aleph_1$-categorical and the upper bound is 1.

The case $d = 2$ only has to be distinguished in the previous subsection. There we were assuming $TT(d-1)$ and $\#RD(d-1)$ is finite so again $T$ is totally transcendental. To obtain an upper bound in this case, we look to subsection 5.6 and so compute only $\#C_1^\alpha$ which in this case (assuming the final property) would be $|\alpha + l|$ for some finite number $l$. So the upper bound becomes

$$|\alpha + \omega|^{|\alpha+1|}$$

which agrees with the general lower bound in this case.

## 6. The spectra

Collecting together all the spectra from Section 4 with the spectra mentioned in the introduction, we obtain the following theorem. This theorem was announced in [8], where examples of theories with each of these spectra were given.

THEOREM 6.1. *For any countable, complete theory $T$ with an infinite model, the uncountable spectrum $\aleph_\alpha \mapsto I(T, \aleph_\alpha)$ ($\alpha > 0$) is the minimum of the map $\aleph_\alpha \mapsto 2^{\aleph_\alpha}$ and one of the following maps:*

1. $2^{\aleph_\alpha}$;
2. $\beth_{d+1}(|\alpha+\omega|)$,           *for some $d$, $\omega \le d < \omega_1$;*
3. $\beth_{d-1}(|\alpha+\omega|^{2^{\aleph_0}})$,     *for some $d$, $0 < d < \omega$;*
4. $\beth_{d-1}(|\alpha+\omega|^{\aleph_0} + \beth_2)$, *for some $d$, $0 < d < \omega$;*
5. $\beth_{d-1}(|\alpha+\omega| + \beth_2)$, *for some $d$, $0 < d < \omega$;*
6. $\beth_{d-1}(|\alpha+\omega|^{\aleph_0})$,      *for some $d$, $0 < d < \omega$;*
7. $\beth_{d-1}(|\alpha+\omega| + 2^{\aleph_0})$, *for some $d$, $1 < d < \omega$;*
8. $\beth_{d-1}(|\alpha+\omega|)$,          *for some $d$, $0 < d < \omega$;*
9. $\beth_{d-2}(|\alpha+\omega|^{|\alpha+1|})$,    *for some $d$, $1 < d < \omega$;*
10. *identically $\beth_2$;*
11. $\begin{cases} |(\alpha+1)^n/\sim_G| - |\alpha^n/\sim_G| & \alpha < \omega, \text{ for some } 1 < n < \omega \text{ and} \\ |\alpha| & \alpha \ge \omega, \text{ some group } G \le \mathrm{Sym}(n); \end{cases}$
12. *identically 1.*

## A. Appendix: On *na*-inclusions

In this first appendix we establish two facts about *na*-inclusions that are used in the text. To simplify the proofs, we extend the definition of an *na*-extension to arbitrary sets. (This more general definition is only used in this appendix.)

*Definition* A.1. If $A \subseteq B$ then $A \subseteq_{na} B$ if whenever $\varphi(x) \in L(A)$ such that $\varphi(B) \setminus A$ is nonempty and $F \subseteq A$ is any finite set then $\varphi(A) \setminus \mathrm{acl}(F)$ is nonempty.

The following lemma follows immediately from a union of chains argument.

LEMMA A.2. *For any set $B$ and any $A \subseteq B$, there is $A' \subseteq_{na} B$ such that $A \subseteq A'$ and $|A'| \le |A| + |L|$.*

LEMMA A.3. *If $B$ is independent from $C$ over $A$ and $A \subseteq_{na} B$ then $C \subseteq_{na} BC$.*

*Proof.* Suppose that we fix $c \in C$ and $b \in B$ so that $\theta(b,c)$ holds and $b \notin \mathrm{acl}(C)$. It suffices to find $b' \in C$ (in fact we will find it in $A$) so that $\theta(b',c)$ holds and $b' \notin \mathrm{acl}(c)$. Let $p = stp(c/B)$. Then $p$ is based on some $a \in A$ and $b$ satisfies $d_p y \theta(x,y)$ that is a formula almost over $a$. Choose $b' \in A \setminus \mathrm{acl}(a)$ which satisfies the same formula. Then immediately $\theta(b',c)$ holds and $b' \notin \mathrm{acl}(c)$. □

LEMMA A.4.  *If $M \subseteq_{na} A$ and $A$ dominates $B$ over $M$ then $M \subseteq_{na} B$.*

*Proof.* Suppose we fix $m \in M$ and $b \in B \setminus M$ so that $\theta(b)$ holds for some $\theta \in L(M)$. Let $c = Cb(stp(b/A)$. $c \in \mathrm{acl}(A) \setminus M$ since $b \notin M$. Moreover, $c \in dcl(\bar{d})$ for a finite sequence of realizations of $\theta$, $\bar{d}$. Choose $c' \in M \setminus \mathrm{acl}(m)$ and $\bar{d}' \in \theta(M)$ so that $c' \in dcl(\bar{d}')$. Since $c' \notin \mathrm{acl}(m)$, one of these realizations must also not be in $\mathrm{acl}(m)$. □

COROLLARY A.5.  *If $M_1$ and $M_2$ are independent over $M_0$, $M_0 \subseteq_{na} M_1$ and $N = M_1 \oplus_{M_0} M_2$ then $M_2 \subseteq_{na} N$.*

## B. Appendix: On $c$- and $d$-isolation

In this appendix we relate two notions of isolation (the first of which is due to Shelah) and give a technical lemma that is used primarily in the proof of Proposition 3.14.

*Definition* B.1 (*T* stable).

- A type $p \in S(A)$ is *c-isolated over $A$* if there is a formula $\varphi(x,b) \in p$ such that no $q \in S(A)$ extending $\varphi(x,b)$ forks over $b$.

- A type $p \in S(A)$ is *d-isolated* (definitionally isolated) *over $A$* if for every $\varphi$ there is a formula $\psi \in p$ such that if for any $q \in S(\mathfrak{C})$ which does not fork over $A$, if $\psi \in q$ then $p|\mathfrak{C}$ and $q$ have the same $\varphi$-definition.

LEMMA B.2 (*T* stable).

1. *If $A$ is algebraically closed and $p \in S(A)$ is c-isolated, then $p$ is d-isolated as well.*

2. *Suppose that $A$ is algebraically closed and $\varphi(a_1,\ldots,a_k)$ holds, where $a_1,\ldots,a_k$ are independent over $A$ and $\mathrm{tp}(a_i/A)$ is d-isolated for each $i$. Then for any sequence $\theta_1,\ldots,\theta_k$ with $\theta_i \in \mathrm{tp}(a_i/A)$ for each $i$, there are $\theta'_i \in \mathrm{tp}(a_i/A)$ that extend $\theta_i$ such that $\varphi(b_1,\ldots,b_k)$ holds for any sequence $b_1,\ldots,b_k$ that are independent over $A$ and satisfy $\theta'_i(b_i)$ for each $i$.*

*Proof.* (1) Let $\psi \in p$ demonstrate that $p$ is c-isolated over $A$ and let $\varphi(x, y)$ be any $L$-formula. Since $A$ is algebraically closed we can choose $\psi \in p$ extending $\psi$ such that $\text{Mult}(\psi', \varphi, \aleph_0) = 1$. Now choose $q \in S(\mathfrak{C})$ that contains $\psi'$ and does not fork over $A$. Since no forking occurs between $\psi$ and $q|A$ we have
$$R(q, \varphi, \aleph_0) = R(q|A, \varphi, \aleph_0) = R(\psi, \varphi, \aleph_0) = R(\psi', \varphi, \aleph_0).$$
Since $\text{Mult}(\psi', \varphi, \aleph_0) = 1$ and $\psi' \in p \cap q$, the $\varphi$-types of $q$ and $p|\mathfrak{C}$ are the same, hence they have the same $\varphi$-definition.

The proof of (2) is by induction on $k$. The case $k = 1$ is clear, so suppose that $k$ is greater than 1. Suppose that $\psi(y_1, \ldots, y_{k-1})$ is the $\varphi$-definition of $tp(a_k/A)$. Since $tp(a_k/A)$ is stationary, $\psi(a_1, \ldots, a_{k-1})$ holds. Now suppose that $\theta_k^*(x)$ is a formula in $tp(a_k/A)$ which is stronger than $\theta_k$ and such that any other type which contains $\theta_k^*$ has the same $\varphi$-definition as $tp(a_k/A)$. Now by induction there are $\theta_1^*, \ldots, \theta_{k-1}^*$ which hold for $a_1, \ldots, a_{k-1}$ respectively, and are such that if $b_1, \ldots, b_{k-1}$ are independent over $A$ and $\theta_i^*(b_i)$ holds for all $i$ then $\psi(b_1, \ldots, b_{k-1})$ holds. Now suppose that $b_1, \ldots, b_k$ are independent over $A$ and $\theta_i^*(b_i)$ holds for all $i$. By the choice of $\theta_k^*$, $\psi$ is the $\varphi$-definition for $tp(b_k/Ab_1, \ldots, b_{k-1})$. By assumption, $\psi(b_1, \ldots, b_{k-1})$ holds and $b_k$ is independent from $b_1, \ldots, b_{k-1}$ over $A$ so $\varphi(b_1, \ldots, b_k)$ holds. $\square$

It is easy to characterize the *c*-isolated types in a superstable theory. To see this, we make the following definition.

*Definition* B.3 (*T* superstable). A consistent $L(A)$-formula $\varphi(x)$ is $R^\infty$-minimal over $A$ if $R^\infty(\varphi) = R^\infty(\varphi')$ for all consistent $L(A)$-formulas $\varphi'(x)$ (with the same free variables) that extend $\varphi$.

If $T$ is superstable then it is readily seen that a type $p \in S(A)$ is c-isolated if and only if it contains a formula that is $R^\infty$-minimal over $A$. The following lemma records three technical results that are used in the text.

LEMMA B.4. *Let $T$ be superstable and let $M \subseteq_{na} \mathfrak{C}$.*

1. *If $M \subseteq A$ and $b$ realizes an $R^\infty$-minimal formula over $A$, then $Ab$ is dominated by $A$ over $M$.*

2. *If $M \subseteq A$ then there is a model $N$ such that $M \subseteq A \subseteq N$ and $N$ is dominated by $A$ over $M$.*

3. *Suppose that $A_1, \ldots, A_k$ are independent over $M$, each $A_i$ is algebraically closed, and $M \subseteq A_i$ for each $i$. If $a_i$ realizes an $R^\infty$-minimal formula $\psi_i$ over $A_i$ for each $i$ and $\varphi(a_1, \ldots, a_k)$, then there are $\psi_i' \in tp(a_i/A_i)$ extending $\psi_i$ for each $i$ such that $\varphi(b_1, \ldots, b_k)$ holds for all $b_1, \ldots, b_k$ such that $\psi_i'(b_i)$ holds for each $i$.*

*Proof.* (1) This is implicit in the proof of Lemma 5.3 of [19].

(2) Since any consistent $L(A)$-formula can be extended to an $R^\infty$-minimal formula over $A$, a model $N$ can be constructed as a sequence of realizations of $c$-isolated types, so the result follows from (1).

(3) Since $\operatorname{tp}(a_i/A_i)$ is $c$-isolated over $A_i$ and $A_i$ is algebraically closed, $\operatorname{tp}(a_i/A_i)$ is $d$-isolated over $A_i$ for each $i$. Further, if $b_1, \ldots, b_k$ are chosen so that $\psi_i(b_i)$ holds for each $i$, then it follows from (1) that $b_i A_i$ is dominated by $A_i$ over $M$. However, since $A_1, \ldots, A_k$ are independent over $M$, this implies that $b_1 A_1, \ldots, b_k A_k$ are independent over $M$ as well. Thus, the result follows from Lemma B.2(2). □

## C. Appendix: $G_\delta$ subsets of Stone spaces

In this appendix we note two facts from descriptive set theory and then establish that several subsets of the Stone space $S(A)$ are $G_\delta$ with respect to the usual topology when $A$ is algebraically closed. These results are used in the proof of Proposition 3.21.

LEMMA C.1. *Let $X$ be any Polish space, i.e., separable, complete metric space.*

1. *Let $E$ be a Borel equivalence relation on $X$. Either $E$ has countably many classes or there is a perfect set of pairwise $E$-inequivalent elements of $X$.*

2. *Any nonempty subset of a countable, $G_\delta$ subset of $X$ has an isolated point.*

*Proof of* 2. Suppose $A = \bigcap_{n \in \omega} U_n$ is nonempty and countable, where each $U_n$ is open in $X$. Since every subset of $A$ is also a $G_\delta$, it suffices to show that $A$ has an isolated point. Let $\bar{A}$ denote the topological closure of $A$ in $X$. As $\bar{A}$ is Polish, it is a Baire space. However, every $U_n \cap \bar{A}$ is dense in $\bar{A}$ and

$$\bigcap_{n \in \omega}(U_n \cap \bar{A}) \cap \bigcap_{a \in A} \bar{A} \smallsetminus \{a\} = \emptyset,$$

so it follows that $\bar{A} \smallsetminus \{a\}$ is not dense in $\bar{A}$ for some $a \in A$. That is, this $a$ is isolated in $A$. □

Fix $A$ countable and algebraically closed. Then the space of types $S(A)$ (with the usual topology) is a Polish space.

LEMMA C.2. $\{p \in S(A) : p \text{ is trivial and } wt(p) = 1\}$ *is a $G_\delta$ subset of* $S(A)$.

*Proof.* We use the following characterization, which is implicit in Theorem 1.8 of [3].

A type $p$ is trivial of weight 1 if and only if

1. For all countable $B \supseteq A$, all weight 1 types $q$ and $r$ over $B$, and all realizations $\langle b, c \rangle$ of $q \otimes r$, if $a$ realizes $p|B$, $a \underset{B}{\perp} b$, and $a \underset{B}{\perp} c$, then $a \underset{B}{\perp} bc$;
2. $p$ has weight at least 1, i.e., for all $B \supseteq A$, for all types $q$, $r$ over $B$, and all realizations $\langle b, c \rangle$ of $q \otimes r$, if $a$ realizes $p|B$ and $a \underset{B}{\not\perp} b$, then $a \underset{B}{\perp} c$; and
3. $p$ is not algebraic.

Now, for all three of these conditions, a specific instance of its negation is witnessed by a formula. Hence, a type $p$ is not trivial, weight 1 if and only if it is not contained in one of the 'bad' formulas. At first glance, it appears like there are continuum constraints, but since the language is countable there really are only countably many. Hence the negation of our set is an $F_\sigma$, so we finish. □

LEMMA C.3. $X_B = \{p : p \text{ determines a complete type over } B\}$ *is a* $G_\delta$ *subset of* $S(A)$ *for every countable* $B \supseteq A$. *Hence,* $\{p : p \perp^a_A q\}$ *is* $G_\delta$ *in* $S(A)$ *for every strong type* $q$ *that is based on* $A$.

*Proof.* The second sentence follows immediately from the first. Fix a countable $B \supseteq A$. For every $L(B)$-formula $\delta(x)$ we say that the $L(A)$-formula *decides* $\delta$ if
$$\varphi(x) \wedge \varphi(y) \quad \text{implies} \quad \delta(x) \leftrightarrow \delta(y).$$
Let $D(\delta) = \{p : \text{there is some } \varphi \in p \text{ that decides } \delta\}$. Clearly, each set $D(\delta)$ is open. Thus, it suffices to show that $p \in X_B$ if and only if $p \in D(\delta)$ for all $\delta(x) \in L(B)$. However, if $p \in X_B$, then
$$p(x) \wedge p(y) \quad \text{implies} \quad \text{tp}(x/B) = \text{tp}(y/B),$$
so it follows by compactness that some $\varphi \in p$ decides $\delta$ for every $\delta \in L(B)$.

Conversely, suppose $p \in D(\delta)$ for all $\delta \in L(B)$. Then, given $a, b$ each realizing $p$, $\delta(a) \leftrightarrow \delta(b)$ for all $\delta \in L(B)$, hence $\text{tp}(a/B) = \text{tp}(b/B)$. □

COROLLARY C.4. $\{p : p \text{ is trivial, } wt(p) = 1 \text{ and } p \perp q\}$ *is a* $G_\delta$ *for every such* $q$.

*Proof.* Recall that if $p$ is trivial and $wt(p) = 1$, then $p \perp q$ if and only if $p \perp^a_A q^{(\omega)}$. □

LEMMA C.5. $\{p \in S(A) : p \underset{A}{\perp^a} A'\}$ is a $G_\delta$ subset of $S(A)$ for any algebraically closed $A' \subseteq A$. In particular,

$$\{p : p \text{ is trivial, } wt(p) = 1 \text{ and } p \perp A'\}$$

is a $G_\delta$ for any such $A'$.

*Proof.* Since $p \underset{A}{\perp^a} A'$ implies $p \perp A'$ among trivial, weight 1 types, the second sentence follows from the first. To conclude the first, say that the $L(A)$-formula $\psi(x)$ *decides* $\delta(x, y) \in L(A)$ if

$$\forall x_1 \forall x_2 \forall y \left( \psi(x_1) \wedge \psi(x_2) \wedge y \underset{A'}{\downarrow} A \to [\delta(x_1, y) \leftrightarrow \delta(x_2, y)] \right)$$

and let $D(\delta) = \{p \in S(A) : \text{some } \psi \in p \text{ decides } \delta(x, y)\}$. Clearly, each $D(\delta)$ is an open subset of $S(A)$. Thus, it suffices to prove that $p \underset{A}{\perp^a} A'$ if and only if $p \in D(\delta)$ for all $\delta(x, y) \in L(A)$.

To see this, first assume that $p \underset{A}{\perp^a} A'$. Then $c \underset{A}{\downarrow} d$ for any $c$ realizing $p$ and any $d \underset{A'}{\downarrow} A$. That is,

$$\forall x_1 \forall x_2 \forall y \left( p(x_1) \wedge p(x_2) \wedge y \underset{A'}{\downarrow} A \to [\mathrm{tp}(x_1 y / A) = \mathrm{tp}(x_2 y / A)] \right),$$

so by compactness, for every $\delta(x, y) \in L(A)$ there is a $\psi(x) \in p$ deciding $\delta$.

Conversely, suppose $p \in D(\delta)$ for every $\delta(x, y) \in L(A)$. Take $c$ realizing $p$ and $d$ freely joined from $A$ over $A'$. Choose $c^*$ realizing $p|Ad$. Since $p \in D(\delta)$ for each $\delta(x, y)$, it follows that $\mathrm{tp}(cd/A) = \mathrm{tp}(c^*d/A)$, so $c \underset{A}{\downarrow} d$ as required. □

MCMASTER UNIVERSITY, HAMILTON, ON, CANADA
*E-mail address*: hartb@mcmcaster.ca

HEBREW UNIVERSITY AT JERUSALEM, JERUSALEM, ISRAEL
*E-mail address*: ehud@math.huji.ac.il

UNIVERSITY OF MARYLAND, COLLEGE PARK, MD
*E-mail address*: mcl@math.umd.edu